\documentclass{article}

\usepackage[utf8]{inputenc}
\usepackage{amsmath, amsthm, amssymb}
\usepackage{longtable}
\usepackage{tikz} 
\usepackage{mathtools}
\usepackage[normalem]{ulem}
\usepackage{setspace}

\usepackage{amsfonts,amscd}	
\usepackage{latexsym}		
\usepackage{graphicx,graphics,epsfig}								
\usepackage{graphpap}										
\usepackage{appendix}	  
								
\usepackage{color}
\usepackage{pdflscape}

\usepackage{xcolor}
\definecolor{maroon}{cmyk}{0,0.87,0.68,0.32}

\definecolor{oran}{RGB}{230,120, 80}
\definecolor{brwn}{RGB}{140, 70, 20}
\definecolor{gren}{RGB}{  0,140, 10}

\newcommand{\up}[1]{\textcolor{orange}{\sf{#1}}}

\newcommand{\gh}[1]{\textcolor{purple}{\sf{#1}}}

\def\bT{{\boldsymbol{T}}}

\newcommand{\eps}[1]{{#1}_{\varepsilon}}

\newtheorem{fact}{Fact}
\newtheorem{lemma}{Lemma}
\newtheorem{theorem}{Theorem}
\newtheorem{corollary}{Corollary}

\begin{document}

\title{On the pebbling numbers of Flower, Blanu\v{s}a and Watkins snarks\footnote{Partially supported by CNPq, FAPERJ and CAPES.}} 

\author{
Matheus Adauto
    \thanks{
        Federal University of Rio de Janeiro, Rio de Janeiro, Brazil}
    \thanks{
        \texttt{adauto@cos.ufrj.br}}
\and
Celina de Figueiredo 
    \footnotemark[1]
    \thanks{
        \texttt{celina@cos.ufrj.br}}
\and
Glenn Hurlbert 
    \thanks{
        Virginia Commonwealth University, Richmond, USA, 
        \texttt{ghurlbert@vcu.edu}}
\and
Diana Sasaki 
    \thanks{
        Rio de Janeiro State University, Rio de Janeiro, Brazil,
        \texttt{diana.sasaki@ime.uerj.br}}
}

\date{}

\maketitle

\doublespacing

\begin{abstract}
Graph pebbling is a game played on graphs with pebbles on their vertices. 
A pebbling move removes two pebbles from one vertex and places one pebble on an adjacent vertex. 
The pebbling number $\pi(G)$ is the smallest $t$ so that from any initial configuration of $t$ pebbles it is possible, after a sequence of pebbling moves, to place a pebble on any given target vertex.
In this paper, we provide the first results on the pebbling numbers of snarks.
Until now, only the Petersen graph had its pebbling number correctly established, although attempts had been made for the Flower and Watkins snarks.

\begin{quote}
\end{quote}
\end{abstract}

\newpage

\section{Introduction} 

Graph pebbling is a mathematical game or puzzle that involves moving pebbles along the edges of a connected graph, subject to certain rules. 
The objective of the game is to place a certain number of pebbles on specific vertices of the graph, typically with the aim of reaching a particular configuration of pebbles or minimizing the number of moves required to achieve a given configuration. 
Various forms of graph pebbling have applications in number theory, computer science, physics, and combinatorial optimization, and have been studied extensively in mathematics (see \cite{HurlbertKenter}). 

In this paper, $G = (V, E)$ is always a simple connected graph. 
The numbers of vertices and edges of $G$ and its diameter are denoted by $n(G)$, $e(G)$, and $D(G)$, respectively.
For a vertex $w$ and positive integer $k$, denote by $N_k[w]$ the set of all vertices that are at a distance at most $k$ from $w$.
The {\it girth} is the length of a shortest cycle in the graph.

\subsection{Pebbling number}
A {\it configuration} $C$ on a graph $G$ is a function $C: V(G) \rightarrow \mathbb{N}$.
The value $C(v)$ represents the number of pebbles at vertex $v$.
A vertex with zero, one, at most one, or at least two pebbles on it is called {\it empty}, {\it a singleton}, {\it small}, or {\it big}, respectively.
The {\it size} $|C|$ of a configuration $C$ is the total number of pebbles on $G$.
A {\it pebbling move} consists of removing two pebbles from a vertex and placing one pebble on an adjacent vertex.
For a {\it target} vertex $r$, $C$ is $r$-{\it solvable} if one can place a pebble on $r$ after a sequence of pebbling moves, and is $r$-{\it unsolvable} otherwise.
It was shown in \cite{HurlKier,MilaClar} that deciding if $C$ is $r$-solvable on $G$ is NP-complete.
The {\it pebbling number} $\pi(G,r)$ is the minimum number $t$ such that every configuration of size $t$ is $r$-solvable.
The {\it pebbling number of} $G$ equals $\pi(G)=\max_r\pi(G,r)$.

The basic lower and upper bounds for every graph are
$\max\{n(G),$ $2^{D(G)}\} \leq \pi(G) \leq (n(G)-D(G))(2^{D(G)}-1)+1$ \cite{ChanGod,glenn1}.
A graph is called \emph{Class 0} if $\pi(G)=n(G)$.
Complete graphs, cubes, the Petersen graph, and many other graphs are known to be Class~0, whereas the cycle graphs $C_n$ satisfy:
$\pi(C_{2d})=2^d$ and $\pi(C_{2d+1})= \lceil (2^{d+2}-1)/3 \rceil$.
It is not yet known whether or not there exist necessary and sufficient conditions for a graph to be Class 0.

\subsection{Snarks}
We define the important family of {\it snark} graphs which are cubic, bridgeless, $4$-edge-chromatic graphs.
They are important for being related to the Four Color Theorem, which holds if and only if no snark is planar~\cite{tait}.  
In \cite{glenn6} we find the origins of the study of the pebbling numbers of chordal graphs.
Here we begin the systematic study of the pebbling numbers of snarks.

The Petersen graph is the smallest snark, having 10 vertices, and was discovered in 1898 \cite{petersen}.
Since then, many others have been discovered.
There are no snarks of order 12, 14, 16, whereas snarks exist for any even order greater than 16.
The Blanu\v{s}a snarks are the two snarks discovered by Danilo Blanu\v{s}a in 1946~\cite{Cavicchioli1998}, when only the Petersen snark was known~\cite{Cavicchioli1998}.
Both Blanu\v{s}a~1 and Blanu\v{s}a~2 snarks have 18 vertices, diameter 4, and girth 5. In this work, they will be denoted by $B_1$ and $B_2$, respectively.
Myriam Preissmann proved in 1982 that there are exactly two snarks of order 18~\cite{Preissmann1982}. We depict in Figure~\ref{fig:Thm_4_example} the Blanu\v{s}a~2 snark together with a vertex labelling that makes clear that the graph was obtained from two copies of the Petersen snark.
See~\cite{Cavicchioli1998} for a thorough history and Table~\ref{snarktable} for a list of several well known snarks.

\begin{figure}
\centering
\begin{tikzpicture}[scale=0.45]
\tikzstyle{every node}=[draw,shape=circle,minimum size=.8cm,font=\scriptsize];
\def \Ang {360/5}
\def \ang {\Ang/4}
\def \radv {2cm}
\def \radz {4.2cm}
\def \radx {6cm}
\path ({90 + (\Ang * 2)}:\radv) node (v2) {$v_2$};
\path ({90 + (\Ang * 1)}:\radv) node (v1) {$v_1$};
\path ({90 + (\Ang * 0)}:\radv) node (v0) {\color{red} $v_0$};
\path ({90 - (\Ang * 1)}:\radv) node (vm1) {$v_{-1}$};
\path ({90 - (\Ang * 2)}:\radv) node (vm2) {$v_{-2}$};
\draw (v0) -- (vm1) -- (vm2) -- (v2) -- (v1) -- (v0);
\path ({90 + (\Ang * 2)}:\radz) node (z2) {$z_2$};
\path ({90 + (\Ang * 1)}:\radz) node (z1) {$z_1$};
\path ({90 + (\Ang * 0)}:\radz) node (z0) {$z_0$};
\node[label={[label distance=-2mm]-90:{\color{gren} $1$}}] at ({90 - (\Ang * 1)}:\radz) (zm1) {$z_{-1}$};
\node[label={[label distance=-2mm]85:{\color{gren} $1$}}] at ({90 - (\Ang * 2)}:\radz) (zm2) {$z_{-2}$};
\draw (v2) -- (z2);
\draw (v1) -- (z1);
\draw (v0) -- (z0);
\draw (vm1) -- (zm1);
\draw (vm2) -- (zm2);
\node[label={[label distance=-2mm]45:{\color{gren} $15$}}] at ({90 + \ang + (\Ang * 2)}:\radx) (x2) {$x_2$};
\path ({90 + \ang + (\Ang * 1)}:\radx) node (x1) {$x_1$};
\path ({90 + \ang + (\Ang * 0)}:\radx) node (x0) {$x_0$};
\node[label={[label distance=-2mm]170:{\color{gren} $1$}}] at ({90 + \ang - (\Ang * 1)}:\radx) (xm1) {$x_{-1}$};
\node[label={[label distance=-2mm]90:{\color{gren} $1$}}] at ({90 + \ang - (\Ang * 2)}:\radx) (xm2) {$x_{-2}$};
\node[label={[label distance=-2mm]90:{\color{gren} $1$}}] at ({90 - \ang + (\Ang * 2)}:\radx) (y2) {$y_2$};
\node[label={[label distance=-2mm]5:{\color{gren} $1$}}] at ({90 - \ang + (\Ang * 1)}:\radx) (y1) {$y_1$};
\path ({90 - \ang + (\Ang * 0)}:\radx) node (y0) {$y_0$};
\node[label={[label distance=-2mm]-100:{\color{gren} $1$}}] at ({90 - \ang - (\Ang * 1)}:\radx) (ym1) {$y_{-1}$};
\path ({90 - \ang - (\Ang * 2)}:\radx) node (ym2) {$y_{-2}$};
\draw (x2) -- (z2) -- (y2);
\draw (x1) -- (z1) -- (y1);
\draw (x0) -- (z0) -- (y0);
\draw (xm1) -- (zm1) -- (ym1);
\draw (xm2) -- (zm2) -- (ym2);
\draw (x2) to[out=-138,in=-150] (x1);
\draw (x1) to[out=150,in=138] (x0);
\draw (x0) to[out=78,in=66] (xm1);
\draw (xm1) to[out=6,in=-6] (xm2);
\draw (xm2) to[out=-90,in=-90] (y2);
\draw (y2) to[out=-174,in=174] (y1);
\draw (y1) to[out=114,in=102] (y0);
\draw (y0) to[out=42,in=30] (ym1);
\draw (ym1) to[out=-30,in=-42] (ym2);
\draw (ym2) to[out=-150,in=-30] (x2);
\end{tikzpicture}
\caption{The graph $J_{5}$ and its (green) $v_0$-unsolvable configuration $C$ of size $22$, which equals the configuration $C^*$ with an extra pebble on $z_{-1}$.}
\label{f:J5}
\end{figure}
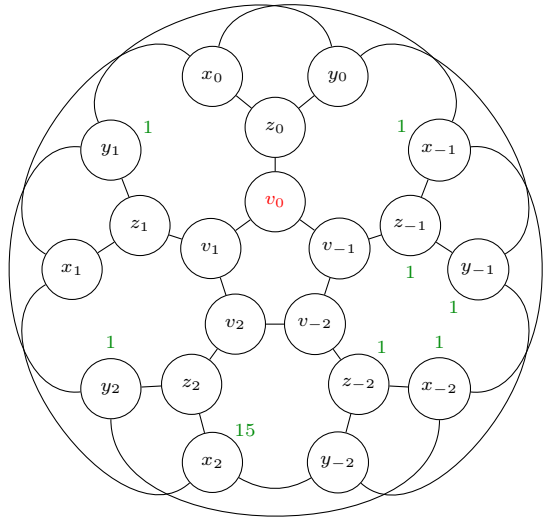

For odd $m=2k+1 \ge 3$, we define the $m^{\rm th}$ {\it Flower snark} $J_m$ as follows (see Figure~\ref{f:J5} for an example with $m=5$)~\cite{Campos}.
For each $i\in\pm\{0,1,\ldots,k\}$ we have vertices $v_i$, $x_i$, and $y_i$ all adjacent to $z_i$.
Thus the number of vertices of the $m^{\rm th}$ Flower snark is $n(J_m)=4m$.
The vertices $\{v_i\}$ form the cycle $C_m$, with adjacencies given by consecutive indices modulo $m$.
The vertices $\{x_i\}$ (resp. $\{y_i\}$) form a path given by the cycle without the edge $x_kx_{-k}$ (resp. $y_ky_{-k}$).
Finally, we add the edges $x_ky_{-k}$ and $y_kx_{-k}$ so that the two paths now form one cycle $C_{2m}$.
It is easy to see that $J_m$ has $m$-fold rotational symmetry, with a necessary twist, along with the reflective symmetry that negates subscripts; that is, the automorphisms of $J_m$ yield three vertex orbits.
Thus the only targets necessary to contemplate are, without loss of generality, $v_0$, $z_0$, and $x_0$.
Thus the Flower snark $J_m$ ($m=2k+1$) has diameter $k+2$, as well as girth $m$ for $m\in\{3,5\}$ and 6 if $m \geq 7$.

The Watkins snark is the snark with 50 vertices discovered by John J. Watkins in 1989~\cite{Watkins1989} depicted in Figure~\ref{f:WatkinsUsual}.
The Watkins snark has diameter 7 and girth 5.

\subsection{Results}

It is known that the Petersen graph is Class 0 \cite{glenn1}.
It is the smallest snark and was the only one whose pebbling number was known. 
We use the Small Neighborhood Lemma presented in Section 2 to prove that the Petersen graph is the only Class 0 snark with at least $23$ vertices or girth at least $5$.

\begin{theorem}
\label{t:Petersen}
The only \emph{Class 0} snark of girth at least 5 is the Petersen graph.
Moreover, if $G$ is a \emph{Class 0} snark with girth at most 4, then $n(G) \le 22$.
\end{theorem}

We also prove the following bounds on the pebbling numbers of snarks.
Recall that the basic lower and upper bounds for a graph are 
$\max\{n(G),$ $2^{D(G)}\} \leq \pi(G) \leq (n(G)-D(G))(2^{D(G)}-1)+1$.
For the Flower snarks, this means that $12 \leq \pi(J_3) \leq 64$,
$20 \leq \pi(J_5) \leq 241$, and
$32 \leq \pi(J_7) \leq 691$.
Theorem~\ref{t:Petersen} improves the $J_5$ lower bound to $21 \leq \pi(J_5)$.
Theorem~\ref{t:Flower} provides much tighter bounds, and corrects a claim of \cite{pebblingflowersnark} that, for $m \geq 5$, $\pi(J_{m})=4m+1$.
In fact, Theorem \ref{t:Flower} identifies the correct order of magnitude ($2^{k+2}$) of $\pi(J_m)$ as $m$ grows, up to some constant between $1$ and $1.8$.

\begin{theorem}
\label{t:Flower}
We have $\pi(J_3) \leq 13$, $23\le \pi(J_5)\le 30$, $41\le \pi(J_7)\le 61$, and, for all $k\ge 3$ with $m=2k+1$, we have $2^{k+2}+9 \le \pi(J_m)\le \lfloor 2^{k+2}9/5 + 2k - 18/5 \rfloor +1$.
\end{theorem}

For the Blanu\v{s}a snarks the basic bounds give $18 \leq \pi(B_i) \leq 211$.
Theorem~\ref{t:Petersen} improves the $B_i$ lower bound to $19 \leq \pi(B_i)$.
Theorem~\ref{t:Blanusa} provides much tighter bounds for $B_i$.

\begin{theorem}
\label{t:Blanusa}
We have $23\leq \pi(B_i)\leq 34$.
\end{theorem}

For the Watkins snark $W$ the basic bounds give $128 \leq \pi(W) \leq 5462$.
Theorem~\ref{t:Watkins} provides a tighter lower bound,
and corrects a claim of \cite{pebblingwatkinssnark} that $\pi(W)=166$.

\begin{theorem}
\label{t:Watkins}
The pebbling number of the Watkins graph $W$ satisfies $183\le \pi(W)$.
\end{theorem}

\section{Techniques}

\subsection{For lower bounds}
\label{ss:lower}

Given a graph $G$ with vertices $u$ and $v$ such that $N_a[u] \cap N_b[v]=\emptyset$ for some non-negative integers $a$ and $b$.
Define the configuration $C^*=C^*_{u,v}$ by $C^*(v)=2^{a+b+1}-1$, $C^*(x)=0$ for all $x\in (N_a[u]\cup N_b[v])-\{v\}$, and $C^*(x)=1$ otherwise.
The authors of \cite{cranston} proved the following lemma (SNL) to provide a lower bound on $\pi(G)$.

\begin{lemma}  (Small Neighborhood Lemma \cite{cranston})
Let $G$ be a graph and $u, v \in V(G)$ such that $N_a[u] \cap N_b[v]=\emptyset$ for some non-negative integers $a$ and $b$.
Then $C^*$ is $u$-unsolvable.
Consequently, $\pi(G)\ge \pi(G,u)> |C^*|$.
In particular, if $N_a[u] \cap N_b[v]=\emptyset$ and $\left|N_a[u] \cup N_b[v]\right|<2^{a+b+1}$, then $G$ is not Class~0.
\label{l:SNL}
\end{lemma}

One can see how SNL is, in some sense, a sharpening of the basic exponential lower bound.
One immediate consequence we will use here is the following corollary.

\begin{corollary}
\label{c:girth5}
If $G$ is an n-vertex cubic graph, with diameter at least 4 and girth at least 5, then there is an unsolvable configuration of size $15 + (n-14)$.
\end{corollary}

Another consequence we will use here is the following corollary.

\begin{corollary}[\cite{cranston}]
\label{c:cranston}
If $G$ is an n-vertex Class 0 graph with diameter at least~3, then $e(G) \geq$ $\frac{5}{3} n-\frac{11}{3}$.
\end{corollary}

A graph $H$ is a {\it retract} of a graph $G$ if there is a function $\phi:V(G)\rightarrow V(H)$ that preserves edges; that is, if $uv\in E(G)$ then $\phi(u)\phi(v)\in E(H)$.
For example, the 5-cycle is a retract of the Petersen graph.
We will make use of the following lemma in the proof 
of the lower bound of Theorem~\ref{t:Blanusa}, since a portion of the Blanu\v{s}a graph has a $C_{9}$ retract,
and
of the lower bound of Theorem~\ref{t:Watkins}, since a portion of the Watkins graph has a $C_{15}$ retract.

\begin{lemma}[Retract Lemma \cite{Chung}]
\label{l:retract}
If $H$ is a retract of $G$, then $\pi(H)\le \pi(G)$.
\end{lemma}

The idea of the proof of Lemma \ref{l:retract} is straightforward, as any solution along edges in $G$ has a corresponding solution along its retracted edges in $H$. 

Another helpful lemma uses the notion of pebble weights.
Given a target vertex $r$, define the $r$-{\it weight} of a pebble on a vertex $v$ to be $2^{-k}$, where $k$ is the distance from $v$ to $r$.
Furthermore, the $r$-{\it weight} of a configuration $C$ is defined to be the sum of the $r$-weights of its pebbles.

\begin{lemma}[\cite{HurlGeneral}]
\label{l:weights}
Suppose that a configuration $C$ on a graph $G$ is $r$-solvable for some target vertex $r$.
Then $C$ has $r$-weight at least 1.
\end{lemma}

The lemma is proved by noting that a pebbling step never increases the $r$-weight of a configuration, and any configuration with a pebble on the target has $r$-weight at least 1.
In fact, the $r$-weight is maintained if and only if the pebbling step decreases the distance to the target $r$.
Hence, when $G$ is a path and $r$ is one of its leaves, $r$-weight at least one characterizes $r$-solvable configurations, which is not true for general graphs.
However, a generalization of this $r$-weight concept is introduced as weights for upper bounds in Section \ref{ss:upper}.

\subsection{For upper bounds}
\label{ss:upper}

Here we describe a linear optimization approach introduced in \cite{HurlbertWFL}.
For an unknown configuration $C$ on a graph $G$ with target vertex $r$, 
we consider a connected subgraph $H$ of $G$ that contains $r$.
The intention is to derive a linear inequality in the variables $C(v)$ with $v\in V(H)$ that is satisfied whenever $C$ is $r$-unsolvable.
Given a collection of such inequalities over various choices of $H$, we can then maximize $|C| = \sum_{v\in V(G)} C(v)$ subject to those constraints, assuming that $C(v)\ge 0$ for all $v\in V(G)$.
The optimum value of this linear program is therefore a strict lower bound on $\pi(G,r)$.
This value may be tight for some graphs; however, this is really an integer optimization problem, and so typically will yield a result less than the truth.
This idea was successfully carried out when $H$ is a tree (generalized to some non-trees in \cite{cranston}).
We introduce the method now.

Let $T$ be a subtree of a graph $G$ rooted at vertex $r$, with at least two vertices. 
For a vertex $v \in V(T)$ let $v^{+}$ denote the {\it parent} of $v$; i.e. the $T$-neighbor of $v$ that is one step closer to $r$ (we also say that $v$ is a {\it child} of $v^{+}$). 
We call $T$ an $r$-{\it strategy} when we associate with it a non-negative {\it weight function} $w_T$ with the property that $w_T(r)=0$ and $w_T\left(v^{+}\right) \ge 2 w_T(v)$ for every other vertex v that is not a neighbor of $r$ (and $w_T(v)=0$ for vertices not in $T$).
Let $\bT$ be the configuration with $\bT(r)=0, \bT(v)=1$ for all other $v \in V(T)$, and $\bT(v)=0$ everywhere else. 
We now define the $T$-{\it weight} of any configuration $C$ (including $\bT$) by $w_T(C)=\sum_{v \in V} w_T(v) C(v)$.
The following lemma (WFL) is used to provide an upper bound on $\pi(G)$.

\begin{lemma}[Weight Function Lemma \cite{HurlbertWFL}]
\label{l:WFL}
Let $T$ be an $r$-strategy of $G$ with associated weight function $w_T$. 
Suppose that $C$ is an $r$-unsolvable configuration of pebbles on $V(G)$. 
Then $w_T(C) \leq w_T(\bT)$.
\end{lemma}

One way to view $T$-weights as a generalization of $r$-weights is as follows.
Consider when $T$ is a path $v_0v_1\cdots v_n$ and let $C$ be a $v_0$-unsolvable configuration on $T$.
Then notice that the formula for $v_0$-weights of pebbles on the vertices of $T$ form a valid weight function $w_T$ if we re-weight $v_0$ at 0 instead of 1, which changes nothing in practice because $C(v_0)=0$.
In this case, WFL states that $\sum_{i=1}^n C(v_i)2^{-i}\le \sum_{i=1}^n 2^{-i} = 1-2^{-n}$.
Because $C$ is a non-negative integer-valued function, this is equivalent to $\sum_{i=1}^n C(v_i)2^{-i} < 1$; i.e. the $v_0$-weight of $C$ is less than one.
Furthermore, the $r$-weights of a general graph $G$ can be thought of as the $T$-weights of a breadth-first-search spanning tree $T$ of $G$, rooted at the target vertex $r$.
Then an $r$-unsolvable configuration $C$ satisfies the WFL inequality and so, by retracting $T$ onto a path of length equal to the eccentricity of $r$, we obtain the aforementioned $r$-weight condition for $C$.

\section{Proofs}

\subsection{Proof of Theorem \ref{t:Petersen}}
\label{s:Petersen}

Note that every snark has exactly $3n/2$ edges, and that the Petersen graph has diameter equal to 2 and all other snarks have diameter at least 3. 
Then, by Corollary \ref{c:cranston}, if $n(G)>22$ we get $3n/2 < (5n-11)/3$.
Therefore, every snark with $n(G)>22$ is not Class $0$.
The remaining non-Petersen snarks with fewer vertices and girth at least $5$ (among them, the Flower $J_5$, the Blanu\v{s}as, and the Loupekines) all have diameter $4>2+1$, so for any vertices $u$ and $v$ at distance 4 from each other we have $|N_2[u]|=10$ and $|N_1[v]|=4$.
Thus $\left|N_2[u] \cup N_1[v]\right| = 14 < 16 = 2^{2+1+1}$, and so none of these graphs are Class $0$ by SNL.
\hfill $\Box$

\subsection{Proof of Theorem \ref{t:Flower}}
\label{s:Flower}

First we prove the lower bounds.
For these we only need to  display a configuration, of size one less than the lower bound, that cannot reach some target.

The diameter of the Flower graph $J_m$, $m=2k+1$, is $k+2$, and the distance between $v_0$ and $x_k$ is $k+2$. The rotational symmetry  of the Flower graph $J_m$ gives $2m$ petals, 
each of which is in a 6-cycle.

For $J_5$, the diameter is 4 and the girth is 5. The $v_0$-unsolvable configuration $C^*_{v_0,x_2}$, that is provided by SNL for $a=2$ and $b=1$ places 15 pebbles on $x_2$ and one pebble on each of the 6 vertices not in set $N_2[v_0] \cap N[x_2]$, and has size 21.
Notice that we can add a pebble to $z_{-1}$ to obtain the configuration $C$ in Figure~\ref{f:J5}.
It is not difficult to argue that $C$ is also $v_0$-unsolvable, since any supposed solution would need to use the pebble at $z_{-1}$.

For $m\ge 7$ (i.e. $k\ge 3$), we will use $C^*=C^*_{v_0,x_k}$ only.
In this case, the girth is equal to 6. 
So $|N_2[v_0]|= 1 + 3 + 6$, and for any integer $3\le i\le k$, the set of vertices at distance $i$ from $v_0$ is $\{v_i, v_{-i}, z_{i+1}, z_{-(i+1)}, x_{i-2}, y_{i-2}, x_{-(i-2)}, y_{-(i-2)} \}$.
Hence, for any integer $2\le i\le k$, we have
$|N_i[v_0]|= 1 + 3 + 6 + (i -2) 8 = 8i-6$.
From this we can compute $|C^*|=(2^{a+b+1}-1) + (n - |N_a[v_0]| - |N_b[x_k]|)$.
In the case that $k$ is even we use $a=k/2$ and $b=a+1$, while if $k$ is odd we use $a=b=(k+1)/2$.
In either case we obtain $|C^*|=2^{k+2}+8$.

Now we prove the upper bounds, using WFL.
We shall define strategies with root $z_0$.
We refer to the Appendix for the strategies with roots $v_0$ and $x_0$, 
that give values not larger than the ones given by strategies presented for root $z_0$.

For $J_3$, we define three $z_0$-strategies $\bT_0$, $\bT_1$, and $\bT_{-1}$ by 
\vspace{-0.04in}
\begin{quote}
    $\bullet$ 
    $\bT_0(v_0,v_1,v_{-1},z_1,z_{-1},x_1,y_1,x_{-1},y_{-1}) = (8,4,4,2,2,1,1,1,1)$, \\
    $\bullet$
    $\bT_1(x_0,x_1,x_{-1},z_1,z_{-1},v_{-1}) = (8,4,4,1,2,1)$ and\\ 
    $\bullet$
    $\bT_{-1}(y_0,y_1, y_{-1},z_1,z_{-1},v_1) = (8,4,4,2,1,1)$,
\end{quote}
giving rise to the inequality $|C| \le \frac{1}{5}(\bT_0+\bT_1+\bT_{-1}) \le 64/5$ whenever $C$ is $v_0$-unsolvable.
Hence $\pi(J_3,z_0)\le 13$.
Using the strategies presented in the Appendix for targets $v_0$ and $x_0$, we may conclude that $\pi(J_3)\le 13$.

For $J_5$, we define three $z_0$-strategies $\bT_0$, $\bT_1$, and $\bT_{-1}$ by 
\vspace{-0.04in}
\begin{quote}
    $\bullet$ 
    $\bT_0(v_0,v_1,v_{-1},v_2,v_{-2},z_1,z_{-1},z_2,z_{-2},x_2,y_2,x_{-2},y_{-2}) = (16,8,8,4,4,4,4,$ $2,2,1,1,1,1)$, \\
    $\bullet$
    $\bT_1(x_0,x_1,x_{-1}, x_2,x_{-2},z_2,z_{-2},v_2,z_1) = (16,8,8,4,4,2,1,1,1)$ and\\ 
    $\bullet$
    $\bT_{-1}(y_0,y_1,y_{-1}, y_2,y_{-2},z_2,z_{-2},v_{-2},z_{-1}) = (16,8,8,4,4,2,1,1,1)$,
\end{quote}
giving rise to the inequality $|C| \le \frac{1}{5}(\bT_0+\bT_1+\bT_{-1}) \le 146/5$ whenever $C$ is $v_0$-unsolvable.
Hence $\pi(J_5,z_0)\le 30$.
Using the strategies presented in the Appendix for targets $v_0$ and $x_0$, we may conclude that $\pi(J_5)\le 30$.

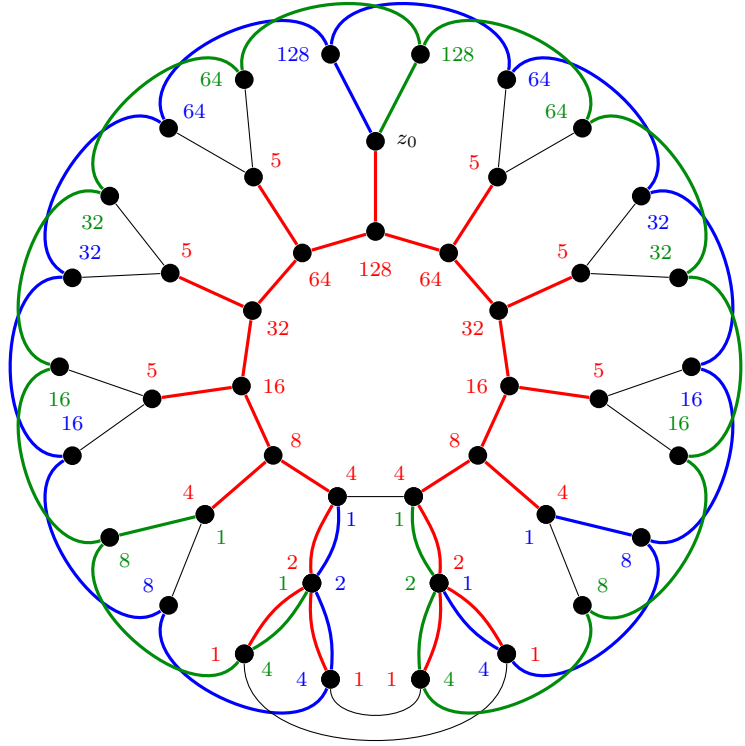
\begin{figure}
\centering
\begin{tikzpicture}[scale=1.2]
\tikzstyle{every node}=[circle,fill=black,inner sep=2.5pt];
\def \Ang {360/11}
\def \ang {\Ang/4}
\def \radv {1.5cm}
\def \radz {2.5cm}
\def \radx {3.5cm}
\node[label={85:{\color{red}\footnotesize $4$}}] at ({90 + (\Ang * 5)}:\radv) (v5) {};
\node[label={-85:{\color{blue}\footnotesize $1$}}] at ({90 + (\Ang * 5)}:\radv) (v5b) {};
\node[label={10:{\color{red}\footnotesize $8$}}] at ({90 + (\Ang * 4)}:\radv) (v4) {};
\node[label={0:{\color{red}\footnotesize $16$}}] at ({90 + (\Ang * 3)}:\radv) (v3) {};
\node[label={-10:{\color{red}\footnotesize $32$}}] at ({90 + (\Ang * 2)}:\radv) (v2) {};
\node[label={280:{\color{red}\footnotesize $64$}}] at ({90 + (\Ang * 1)}:\radv) (v1) {};
\node[label={270:{\color{red}\footnotesize $128$}}] at ({90 + (\Ang * 0)}:\radv) (v0) {};
\node[label={260:{\color{red}\footnotesize $64$}}] at ({90 - (\Ang * 1)}:\radv) (vm1) {};
\node[label={190:{\color{red}\footnotesize $32$}}] at ({90 - (\Ang * 2)}:\radv) (vm2) {};
\node[label={180:{\color{red}\footnotesize $16$}}] at ({90 - (\Ang * 3)}:\radv) (vm3) {};
\node[label={170:{\color{red}\footnotesize $8$}}] at ({90 - (\Ang * 4)}:\radv) (vm4) {};
\node[label={95:{\color{red}\footnotesize $4$}}] at ({90 - (\Ang * 5)}:\radv) (vm5) {};
\node[label={-95:{\color{gren}\footnotesize $1$}}] at ({90 - (\Ang * 5)}:\radv) (vm5g) {};
\draw[very thick,red] (v0) -- (vm1) -- (vm2) -- (vm3) -- (vm4) -- (vm5);
\draw[very thick,red] (v5) -- (v4) -- (v3) -- (v2) -- (v1) -- (v0);
\draw (vm5) -- (v5);
\node[label={0:{\color{blue}\footnotesize $2$}}] at ({90 + (\Ang * 5)}:\radz) (z5) {};
\node[label={130:{\color{red}\footnotesize $2$}}] at ({90 + (\Ang * 5)}:\radz) (z5r) {};
\node[label={180:{\color{gren}\footnotesize $1$}}] at ({90 + (\Ang * 5)}:\radz) (z5g) {};
\node[label={110:{\color{red}\footnotesize $4$}}] at ({90 + (\Ang * 4)}:\radz) (z4) {};
\node[label={290:{\color{gren}\footnotesize $1$}}] at ({90 + (\Ang * 4)}:\radz) (z4g) {};
\node[label={90:{\color{red}\footnotesize $5$}}] at ({90 + (\Ang * 3)}:\radz) (z3) {};
\node[label={70:{\color{red}\footnotesize $5$}}] at ({90 + (\Ang * 2)}:\radz) (z2) {};
\node[label={20:{\color{red}\footnotesize $5$}}] at ({90 + (\Ang * 1)}:\radz) (z1) {};
\node[label={0:{\footnotesize $z_0$}}] at ({90 - (\Ang * 0)}:\radz) (z0) {};
\node[label={170:{\color{red}\footnotesize $5$}}] at ({90 - (\Ang * 1)}:\radz) (zm1) {};
\node[label={120:{\color{red}\footnotesize $5$}}] at ({90 - (\Ang * 2)}:\radz) (zm2) {};
\node[label={90:{\color{red}\footnotesize $5$}}] at ({90 - (\Ang * 3)}:\radz) (zm3) {};
\node[label={70:{\color{red}\footnotesize $4$}}] at ({90 - (\Ang * 4)}:\radz) (zm4) {};
\node[label={250:{\color{blue}\footnotesize $1$}}] at ({90 - (\Ang * 4)}:\radz) (zm4b) {};
\node[label={180:{\color{gren}\footnotesize $2$}}] at ({90 - (\Ang * 5)}:\radz) (zm5) {};
\node[label={50:{\color{red}\footnotesize $2$}}] at ({90 - (\Ang * 5)}:\radz) (zm5g) {};
\node[label={0:{\color{blue}\footnotesize $1$}}] at ({90 - (\Ang * 5)}:\radz) (zm5b) {};
\draw[very thick,blue] (v5) to[out=275,in=55] (z5);
\draw[very thick,red] (v5) to[out=240,in=95] (z5);
\draw[very thick,red] (v4) -- (z4);
\draw[very thick,red] (v3) -- (z3);
\draw[very thick,red] (v2) -- (z2);
\draw[very thick,red] (v1) -- (z1);
\draw[very thick,red] (v0) -- (z0);
\draw[very thick,red] (vm1) -- (zm1);
\draw[very thick,red] (vm2) -- (zm2);
\draw[very thick,red] (vm3) -- (zm3);
\draw[very thick,red] (vm4) -- (zm4);
\draw[very thick,gren] (vm5) to[out=265,in=125] (zm5);
\draw[very thick,red] (vm5) to[out=300,in=85] (zm5);
\node[label={180:{\color{blue}\footnotesize $4$}}] at ({90 + \ang + (\Ang * 5)}:\radx) (x5) {};
\node[label={0:{\color{red}\footnotesize $1$}}] at ({90 + \ang + (\Ang * 5)}:\radx) (x5r) {};
\node[label={135:{\color{blue}\footnotesize $8$}}] at ({90 + \ang + (\Ang * 4)}:\radx) (x4) {};
\node[label={90:{\color{blue}\footnotesize $16$}}] at ({90 + \ang + (\Ang * 3)}:\radx) (x3) {};
\node[label={80:{\color{blue}\footnotesize $32$}}] at ({90 + \ang + (\Ang * 2)}:\radx) (x2) {};
\node[label={5:{\color{blue}\footnotesize $64$}}] at ({90 + \ang + (\Ang * 1)}:\radx) (x1) {};
\node[label={180:{\color{blue}\footnotesize $128$}}] at ({90 + \ang + (\Ang * 0)}:\radx) (x0) {};
\node[label={0:{\color{blue}\footnotesize $64$}}] at ({90 + \ang - (\Ang * 1)}:\radx) (xm1) {};
\node[label={-85:{\color{blue}\footnotesize $32$}}] at ({90 + \ang - (\Ang * 2)}:\radx) (xm2) {};
\node[label={-90:{\color{blue}\footnotesize $16$}}] at ({90 + \ang - (\Ang * 3)}:\radx) (xm3) {};
\node[label={-100:{\color{blue}\footnotesize $8$}}] at ({90 + \ang - (\Ang * 4)}:\radx) (xm4) {};
\node[label={190:{\color{blue}\footnotesize $4$}}] at ({90 + \ang - (\Ang * 5)}:\radx) (xm5) {};
\node[label={0:{\color{red}\footnotesize $1$}}] at ({90 + \ang - (\Ang * 5)}:\radx) (xm5r) {};
\node[label={-10:{\color{gren}\footnotesize $4$}}] at ({90 - \ang + (\Ang * 5)}:\radx) (y5) {};
\node[label={180:{\color{red}\footnotesize $1$}}] at ({90 - \ang + (\Ang * 5)}:\radx) (y5r) {};
\node[label={-80:{\color{gren}\footnotesize $8$}}] at ({90 - \ang + (\Ang * 4)}:\radx) (y4) {};
\node[label={-90:{\color{gren}\footnotesize $16$}}] at ({90 - \ang + (\Ang * 3)}:\radx) (y3) {};
\node[label={265:{\color{gren}\footnotesize $32$}}] at ({90 - \ang + (\Ang * 2)}:\radx) (y2) {};
\node[label={180:{\color{gren}\footnotesize $64$}}] at ({90 - \ang + (\Ang * 1)}:\radx) (y1) {};
\node[label={0:{\color{gren}\footnotesize $128$}}] at ({90 - \ang - (\Ang * 0)}:\radx) (y0) {};
\node[label={175:{\color{gren}\footnotesize $64$}}] at ({90 - \ang - (\Ang * 1)}:\radx) (ym1) {};
\node[label={100:{\color{gren}\footnotesize $32$}}] at ({90 - \ang - (\Ang * 2)}:\radx) (ym2) {};
\node[label={90:{\color{gren}\footnotesize $16$}}] at ({90 - \ang - (\Ang * 3)}:\radx) (ym3) {};
\node[label={45:{\color{gren}\footnotesize $8$}}] at ({90 - \ang - (\Ang * 4)}:\radx) (ym4) {};
\node[label={0:{\color{gren}\footnotesize $4$}}] at ({90 - \ang - (\Ang * 5)}:\radx) (ym5) {};
\node[label={180:{\color{red}\footnotesize $1$}}] at ({90 - \ang - (\Ang * 5)}:\radx) (ym5r) {};
\draw[very thick,red] (y5) to[out=60,in=210] (z5);
\draw[very thick,gren] (y5) to[out=30,in=240] (z5);
\draw[very thick,red] (x5) to[out=115,in=265] (z5);
\draw[very thick,blue] (x5) to[out=85,in=295] (z5);
\draw (x4) -- (z4);
\draw[very thick,gren] (z4) -- (y4);
\draw (x3) -- (z3) -- (y3);
\draw (x2) -- (z2) -- (y2);
\draw (x1) -- (z1) -- (y1);
\draw[very thick,blue] (x0) -- (z0);
\draw[very thick,gren] (z0) -- (y0);
\draw (xm1) -- (zm1) -- (ym1);
\draw (xm2) -- (zm2) -- (ym2);
\draw (xm3) -- (zm3) -- (ym3);
\draw[very thick,blue] (xm4) -- (zm4);
\draw (zm4) -- (ym4);
\draw[very thick,red] (zm5) to[out=330,in=120] (xm5);
\draw[very thick,blue] (zm5) to[out=300,in=150] (xm5);
\draw[very thick,red] (zm5) to[out=275,in=65] (ym5);
\draw[very thick,gren] (zm5) to[out=245,in=95] (ym5);
\draw[very thick,blue] (x5) to[out=248,in=243] (x4);
\draw[very thick,blue] (x4) to[out=215,in=210] (x3);
\draw[very thick,blue] (x3) to[out=183,in=177] (x2);
\draw[very thick,blue] (x2) to[out=150,in=145] (x1);
\draw[very thick,blue] (x1) to[out=117,in=112] (x0);
\draw[very thick,blue] (x0) to[out=85,in=79] (xm1);
\draw[very thick,blue] (xm1) to[out=52,in=46] (xm2);
\draw[very thick,blue] (xm2) to[out=19,in=14] (xm3);
\draw[very thick,blue] (xm3) to[out=-14,in=-19] (xm4);
\draw[very thick,blue] (xm4) to[out=-46,in=-52] (xm5);
\draw (xm5) to[out=-90,in=-90] (y5);
\draw[very thick,gren] (y5) to[out=232,in=226] (y4);
\draw[very thick,gren] (y4) to[out=199,in=194] (y3);
\draw[very thick,gren] (y3) to[out=166,in=161] (y2);
\draw[very thick,gren] (y2) to[out=134,in=128] (y1);
\draw[very thick,gren] (y1) to[out=101,in=95] (y0);
\draw[very thick,gren] (y0) to[out=68,in=63] (ym1);
\draw[very thick,gren] (ym1) to[out=35,in=30] (ym2);
\draw[very thick,gren] (ym2) to[out=3,in=-3] (ym3);
\draw[very thick,gren] (ym3) to[out=-30,in=-36] (ym4);
\draw[very thick,gren] (ym4) to[out=-63,in=-68] (ym5);
\draw (ym5) to[out=-90,in=-90] (x5);
\end{tikzpicture}
\caption{The graph $J_{11}$ and its three $z_0$-strategies $\bT_0$ (in red), $\bT_1$ (in blue), and $\bT_{-1}$ (in green).}
\label{f:J11}
\end{figure}

For $m\ge 5$ (i.e. $k\ge 2$), please refer to Figure~\ref{f:J11} to see $m = 11$.
Using the same pattern we have defined above for the three $z_0$-strategies for $J_5$, we define three corresponding $z_0$-strategies by
\begin{itemize}
    \item 
    $\bT_0(v_0,v_1,v_{-1},v_2,v_{-2},\ldots, v_k,v_{-k},z_k,z_{-k},x_k,y_k,x_{-k},y_{-k})$\\
    $= (2^{k+2},2^{k+1},2^{k+1},2^k,2^k,\ldots, 2,2,1,1,1,1)$ and\\
    $\bT_0(z_1,z_{-1},\ldots,z_{k-2},z_{2-k},z_{k-1},z_{1-k})
    = (5,5,\ldots,5,5,4,4)$;
    \item
    $\bT_1(x_0,x_1,x_{-1},\ldots,x_k,x_{-k},z_k,v_k)\\
    = (2^{k+2},2^{k+1},2^{k+1},\ldots,4,4,2,1)$ and
    $\bT_1(z_{2-k},z_{1-k}) = (1,1)$; and
    \item
    $\bT_{-1}(y_0,y_1,y_{-1},\ldots,y_k,y_{-k},z_{-k},v_{-k})\\
    = (2^{k+2},2^{k+1},2^{k+1},\ldots,4,4,2,1)$ and
    $\bT_{-1}(z_{k-2},z_{k-1}) = (1,1)$.
\end{itemize}

The sum $\bT_0+\bT_1+\bT_{-1}$ has $3$ vertices with coefficient $2^{k+2}$, $6$ with $2^i$ (for each $3\le i\le k+1$), and $2k+6$ with coefficient $5$, giving rise to the inequality 
\begin{align*}
    5|C|
    &\le \bT_0+\bT_1+\bT_{-1}\\
    &= 3(2^{k+2}) + 6(2^3 + \cdots + 2^{k+1})  + 5(2k+6)\\
    &= 6(2^3 + \cdots + 2^{k+2}) - 3(2^{k+2}) + 5(2k+6)\\
    &= 48(2^k-1) - 3(2^{k+2}) + 10k + 30\\
    &= 36(2^{k}) + 10k - 18\\
    &= 9(2^{k+2}) + 10k - 18,
\end{align*}
whenever $C$ is $z_0$-unsolvable.
Hence $|C|\le 2^{k+2}9/5 + 2k - 18/5$, and so 
$\pi(J_m,z_0) \le \lfloor 2^{k+2}9/5 + 2k - 18/5 \rfloor +1$.
Using the strategies presented in the Appendix for targets $v_0$ and $x_0$, we may conclude that 
$\pi(J_m)\le \lfloor 2^{k+2}9/5 + 2k - 18/5 \rfloor +1$.
In particular, $\pi(J_7)\le 61$.
\hfill $\Box$

\subsection{Proof of Theorem \ref{t:Blanusa}}

First we prove the lower bounds. 
Corollary~\ref{c:girth5} gives the same lower bound $20 \leq \pi(G)$, for an arbitrary diameter 4, girth 5, cubic graph $G$ with 18 vertices, like  $B_1$ and $B_2$, by defining an unsolvable configuration $C^*$ of size 19, as follows.
Consider a pair of vertices $u$ and $v$ at distance 4 in such a graph $G$. 
Define $C^*(v)=15$, $C^*(x)=0$ for all $x\in (N_2[u]\cup N_1[v])-\{v\}$, and $C^*(x)=1$ otherwise. 
Hence, $|C^*|=19$.

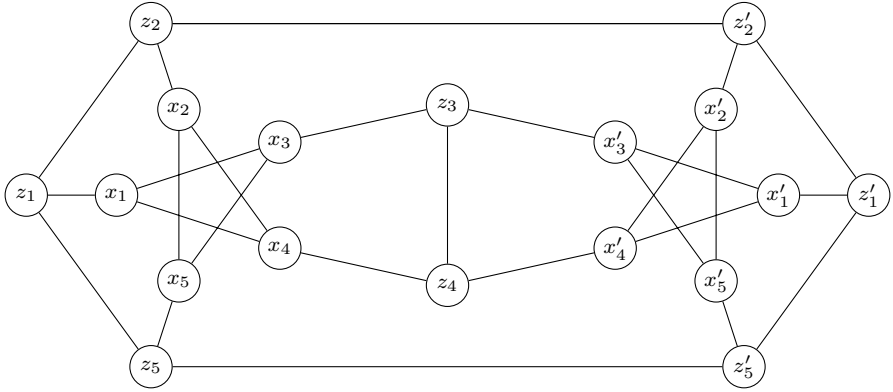
\begin{figure}[ht] 
\centering
\begin{tikzpicture} [scale=0.3]

\def \dia {16}
\def \Ang {360/5}
\def \rada {4}
\def \radb {8}
\def \dist {.4*\radb}
\def \height {2}
\node[draw=none,fill=none](C) at (0:0) {};
\node[draw,circle,minimum size=\dia,inner sep=0,shift=(180:\dist)](x1) at ({180 - (\Ang * 0)}:\rada) {\footnotesize $x_{1}$};
\node[draw,circle,minimum size=\dia,inner sep=0,shift=(180:\dist)](x2) at ({180 - (\Ang * 1)}:\rada) {\footnotesize $x_{2}$};
\node[draw,circle,minimum size=\dia,inner sep=0,shift=(180:\dist)](x3) at ({180 - (\Ang * 2)}:\rada) {\footnotesize $x_{3}$};
\node[draw,circle,minimum size=\dia,inner sep=0,shift=(180:\dist)](x4) at ({180 - (\Ang * 3)}:\rada) {\footnotesize $x_{4}$};
\node[draw,circle,minimum size=\dia,inner sep=0,shift=(180:\dist)](x5) at ({180 - (\Ang * 4)}:\rada) {\footnotesize $x_{5}$};
\node[draw,circle,minimum size=\dia,inner sep=0,shift=(180:\dist)](z1) at ({180 - (\Ang * 0)}:\radb) {\footnotesize $z_{1}$};
\node[draw,circle,minimum size=\dia,inner sep=0,shift=(180:\dist)](z2) at ({180 - (\Ang * 1)}:\radb) {\footnotesize $z_{2}$};
\node[draw,circle,minimum size=\dia,inner sep=0](z3) at ({90 - (\Ang * 0)}:\rada) {\footnotesize $z_{3}$};
\node[draw,circle,minimum size=\dia,inner sep=0](z4) at ({-90 - (\Ang * 0)}:\rada) {\footnotesize $z_{4}$};
\node[draw,circle,minimum size=\dia,inner sep=0,shift=(180:\dist)](z5) at ({180 - (\Ang * 4)}:\radb) {\footnotesize $z_{5}$};
\node[draw,circle,minimum size=\dia,inner sep=0,shift=(0:\dist)](x1p) at ({0 + (\Ang * 0)}:\rada) {\footnotesize $x'_{1}$};
\node[draw,circle,minimum size=\dia,inner sep=0,shift=(0:\dist)](x2p) at ({0 + (\Ang * 1)}:\rada) {\footnotesize $x'_{2}$};
\node[draw,circle,minimum size=\dia,inner sep=0,shift=(0:\dist)](x3p) at ({0 + (\Ang * 2)}:\rada) {\footnotesize $x'_{3}$};
\node[draw,circle,minimum size=\dia,inner sep=0,shift=(0:\dist)](x4p) at ({0 + (\Ang * 3)}:\rada) {\footnotesize $x'_{4}$};
\node[draw,circle,minimum size=\dia,inner sep=0,shift=(0:\dist)](x5p) at ({0 + (\Ang * 4)}:\rada) {\footnotesize $x'_{5}$};
\node[draw,circle,minimum size=\dia,inner sep=0,shift=(0:\dist)](z1p) at ({0 + (\Ang * 0)}:\radb) {\footnotesize $z'_{1}$};
\node[draw,circle,minimum size=\dia,inner sep=0,shift=(0:\dist)](z2p) at ({0 + (\Ang * 1)}:\radb) {\footnotesize $z'_{2}$};
\node[draw,circle,minimum size=\dia,inner sep=0,shift=(0:\dist)](z5p) at ({0 + (\Ang * 4)}:\radb) {\footnotesize $z'_{5}$};

\draw (x1) -- (z1);
\draw (x2) -- (z2);
\draw (x3) -- (z3);
\draw (x4) -- (z4);
\draw (x5) -- (z5);
\draw (x1p) -- (z1p);
\draw (x2p) -- (z2p);
\draw (x3p) -- (z3);
\draw (x4p) -- (z4);
\draw (x5p) -- (z5p);
\draw (x1) -- (x3) -- (x5) -- (x2) -- (x4) -- (x1);
\draw (x1p) -- (x3p) -- (x5p) -- (x2p) -- (x4p) -- (x1p);
\draw (z1) -- (z2);
\draw (z1) -- (z5);
\draw (z1p) -- (z2p);
\draw (z1p) -- (z5p);
\draw (z2) -- (z2p);
\draw (z2) -- (z2p);
\draw (z2) -- (z2p);
\draw (z3) -- (z4);
\draw (z5) -- (z5p);

\end{tikzpicture}
\caption{Blanu\v{s}a 2 and its labelling.} 
\label{fig:Thm_4_example}
\end{figure}

Please refer to Figure~\ref{fig:Thm_4_example}. Analogous arguments can be done for $B_1$ and $B_2$, so we present with no loss of generality the arguments for $B_2$. We can establish the non trivial lower bound $21 \leq \pi(B_2)$, as an application of the Retract Lemma~\ref{l:retract},
since a portion of $B_2$ has a $C_{9}$ retract and $\pi(C_{9}) = 21$.
We can actually define a $x_3$-unsolvable configuration of size 22, by organizing $B_2$ by distance from target $x_3$. 
Consider the $C_{9}$ induced by $x_3, x_1, z_1, z_2, z'_2, x'_2, x'_5, x'_3, z_3$. Adjacent vertices $z'_2$ and $x'_2$ are at distance 4 from target $x_3$. Place 10 pebbles on $z'_2$ and 10 pebbles on $x'_2$ to get a $x_3$-unsolvable configuration of size 20. Additionally, place 1 pebble on $z'_5$ and 1 pebble on $x'_1$, to get the desired $x_3$-unsolvable configuration of size 22, since $z'_5$ and $x'_1$ are among the vertices at distance 3 from target $x_3$, that are not in a cycle $C_{9}$ together with $x_3, z'_2$ and $x'_2$. Hence $23 \leq \pi(B_2)$.

Now we prove the upper bound $\pi(B_2) \leq 34$, using WFL. Note that we have six different roots, by considering: $z_1=z'_1$, $z_2=z'_2=z_5=z'_5$, $x_1=x'_1$, $x_2=x'_2=x_5=x'_5$, $x_3=x'_3=x_4=x'_4$, and $z_3=z_4$.
We shall define strategies with root $x_3$, since this root gave us the largest upper bound. 
We refer to the Appendix for the strategies with the other five roots that give values not larger than the ones given by strategies presented for root $x_3$.

For $B_2$, we define three $x_3$-strategies $\bT_1$, $\bT_2$, and $\bT_3$ by 
\vspace{-0.04in}
\begin{quote}
    $\bullet$ 
    $\bT_1(x_1,z_1,x_4,z_2,z_4,z'_2,x'_4,x'_2) = (32,16,16,8,8,4,4,2)$, \\
    $\bullet$
    $\bT_2(x_5,x_2,z_5,z'_5,z'_1,z'_2,x'_2) = (32,7,16,8,4,2,1)$ and\\ 
    $\bullet$
    $\bT_3(z_3,x'_3, x'_1,x'_5,z'_1,x'_2,x'_4) = (32,16,8,8,4,4,3)$,
\end{quote}
giving rise to the inequality $|C| \le \frac{1}{7}(\bT_1+\bT_2+\bT_3) \le 236/7$ whenever $C$ is $x_3$-unsolvable.
Hence $\pi(B_2,x_3)\le 34$.
Using the strategies presented in the
Appendix for the other five targets, we may conclude that $\pi(B_2) \leq 34$.
\hfill $\Box$

\subsection{Proof of Theorem~\ref{t:Watkins}}
\label{s:Watkins}

\begin{figure}
\centering
\begin{tikzpicture}[scale=0.45]
\def \dia {16}
\def \Ang {360/5}
\def \anga {15}
\def \rada {6}
\def \radb {8}
\def \radc {10}
\def \radd {12}
\def \rade {14}
\node[draw=none,fill=none](U0) at ({90 + (\Ang * 0)}:\rade) {\small $U_0$};
\node[draw,circle,minimum size=\dia,inner sep=0](a1) at ({90 + (\Ang * 0) + \anga}:\radd) {\footnotesize $a_{1}$};
\node[draw,circle,minimum size=\dia,inner sep=0](a2) at ({90 + (\Ang * 0) }:\radd) {\footnotesize $a_{2}$};
\node[draw,circle,minimum size=\dia,inner sep=0](a3) at ({90 + (\Ang * 0) - \anga}:\radd) {\footnotesize $a_{3}$};
\node[draw,circle,minimum size=\dia,inner sep=0](b1) at ({90 + (\Ang * 0) + \anga}:\radc) {\footnotesize $b_{1}$};
\node[draw,circle,minimum size=\dia,inner sep=0](b3) at ({90 + (\Ang * 0) - \anga}:\radc) {\footnotesize $b_{3}$};
\node[draw,circle,minimum size=\dia,inner sep=0](b19) at ({90 + (\Ang * 0) + \anga}:\radb) {\footnotesize $b_{19}$};
\node[draw,circle,minimum size=\dia,inner sep=0](b2) at ({90 + (\Ang * 0) }:\radb) {\footnotesize $b_{2}$};
\node[draw,circle,minimum size=\dia,inner sep=0](b20) at ({90 + (\Ang * 0) - \anga}:\radb) {\footnotesize $b_{20}$};
\node[draw,circle,minimum size=\dia,inner sep=0](a19) at ({90 + (\Ang * 0) + \anga}:\rada) {\footnotesize $a_{19}$};
\node[draw,circle,minimum size=\dia,inner sep=0](a20) at ({90 + (\Ang * 0) - \anga}:\rada) {\footnotesize $a_{20}$};
\node[draw=none,fill=none](U1) at ({90 - (\Ang * 1)}:\rade) {\small $U_1$};
\node[draw,circle,minimum size=\dia,inner sep=0](a11) at ({90 - (\Ang * 1) + \anga}:\radd) {\footnotesize $a_{11}$};
\node[draw,circle,minimum size=\dia,inner sep=0](a12) at ({90 - (\Ang * 1) }:\radd) {\footnotesize $a_{12}$};
\node[draw,circle,minimum size=\dia,inner sep=0](a13) at ({90 - (\Ang * 1) - \anga}:\radd) {\footnotesize $a_{13}$};
\node[draw,circle,minimum size=\dia,inner sep=0](b11) at ({90 - (\Ang * 1) + \anga}:\radc) {\footnotesize $b_{11}$};
\node[draw,circle,minimum size=\dia,inner sep=0](b13) at ({90 - (\Ang * 1) - \anga}:\radc) {\footnotesize $b_{13}$};
\node[draw,circle,minimum size=\dia,inner sep=0](b4) at ({90 - (\Ang * 1) + \anga}:\radb) {\footnotesize $b_{4}$};
\node[draw,circle,minimum size=\dia,inner sep=0](b12) at ({90 - (\Ang * 1) }:\radb) {\footnotesize $b_{12}$};
\node[draw,circle,minimum size=\dia,inner sep=0](b5) at ({90 - (\Ang * 1) - \anga}:\radb) {\footnotesize $b_{5}$};
\node[draw,circle,minimum size=\dia,inner sep=0](a4) at ({90 - (\Ang * 1) + \anga}:\rada) {\footnotesize $a_{4}$};
\node[draw,circle,minimum size=\dia,inner sep=0](a5) at ({90 - (\Ang * 1) - \anga}:\rada) {\footnotesize $a_{5}$};
\node[draw=none,fill=none](U4) at ({90 + (\Ang * 1)}:\rade) {\small $U_4$};
\node[draw,circle,minimum size=\dia,inner sep=0](a16) at ({90 + (\Ang * 1) + \anga}:\radd) {\footnotesize $a_{16}$};
\node[draw,circle,minimum size=\dia,inner sep=0](a17) at ({90 + (\Ang * 1) }:\radd) {\footnotesize $a_{17}$};
\node[draw,circle,minimum size=\dia,inner sep=0](a18) at ({90 + (\Ang * 1) - \anga}:\radd) {\footnotesize $a_{18}$};
\node[draw,circle,minimum size=\dia,inner sep=0](b16) at ({90 + (\Ang * 1) + \anga}:\radc) {\footnotesize $b_{16}$};
\node[draw,circle,minimum size=\dia,inner sep=0](b18) at ({90 + (\Ang * 1) - \anga}:\radc) {\footnotesize $b_{18}$};
\node[draw,circle,minimum size=\dia,inner sep=0](b9) at ({90 + (\Ang * 1) + \anga}:\radb) {\footnotesize $b_{9}$};
\node[draw,circle,minimum size=\dia,inner sep=0](b17) at ({90 + (\Ang * 1) }:\radb) {\footnotesize $b_{17}$};
\node[draw,circle,minimum size=\dia,inner sep=0](b10) at ({90 + (\Ang * 1) - \anga}:\radb) {\footnotesize $b_{10}$};
\node[draw,circle,minimum size=\dia,inner sep=0](a9) at ({90 + (\Ang * 1) + \anga}:\rada) {\footnotesize $a_{9}$};
\node[draw,circle,minimum size=\dia,inner sep=0](a10) at ({90 + (\Ang * 1) - \anga}:\rada) {\footnotesize $a_{10}$};
\node[draw=none,fill=none](U2) at ({90 - (\Ang * 2)}:\rade) {\small $U_2$};
\node[draw,circle,minimum size=\dia,inner sep=0](a21) at ({90 - (\Ang * 2) + \anga}:\radd) {\footnotesize $a_{21}$};
\node[draw,circle,minimum size=\dia,inner sep=0](a22) at ({90 - (\Ang * 2) }:\radd) {\footnotesize $a_{22}$};
\node[draw,circle,minimum size=\dia,inner sep=0](a23) at ({90 - (\Ang * 2) - \anga}:\radd) {\footnotesize $a_{23}$};
\node[draw,circle,minimum size=\dia,inner sep=0](b21) at ({90 - (\Ang * 2) + \anga}:\radc) {\footnotesize $b_{21}$};
\node[draw,circle,minimum size=\dia,inner sep=0](b23) at ({90 - (\Ang * 2) - \anga}:\radc) {\footnotesize $b_{23}$};
\node[draw,circle,minimum size=\dia,inner sep=0](b14) at ({90 - (\Ang * 2) + \anga}:\radb) {\footnotesize $b_{14}$};
\node[draw,circle,minimum size=\dia,inner sep=0](b22) at ({90 - (\Ang * 2) }:\radb) {\footnotesize $b_{22}$};
\node[draw,circle,minimum size=\dia,inner sep=0](b15) at ({90 - (\Ang * 2) - \anga}:\radb) {\footnotesize $b_{15}$};
\node[draw,circle,minimum size=\dia,inner sep=0](a14) at ({90 - (\Ang * 2) + \anga}:\rada) {\footnotesize $a_{14}$};
\node[draw,circle,minimum size=\dia,inner sep=0](a15) at ({90 - (\Ang * 2) - \anga}:\rada) {\footnotesize $a_{15}$};
\node[draw=none,fill=none](U3) at ({90 + (\Ang * 2)}:\rade) {\small $U_3$};
\node[draw,circle,minimum size=\dia,inner sep=0](a6) at ({90 + (\Ang * 2) + \anga}:\radd) {\footnotesize $a_{6}$};
\node[draw,circle,minimum size=\dia,inner sep=0](a7) at ({90 + (\Ang * 2) }:\radd) {\footnotesize $a_{7}$};
\node[draw,circle,minimum size=\dia,inner sep=0](a8) at ({90 + (\Ang * 2) - \anga}:\radd) {\footnotesize $a_{8}$};
\node[draw,circle,minimum size=\dia,inner sep=0](b6) at ({90 + (\Ang * 2) + \anga}:\radc) {\footnotesize $b_{6}$};
\node[draw,circle,minimum size=\dia,inner sep=0](b8) at ({90 + (\Ang * 2) - \anga}:\radc) {\footnotesize $b_{8}$};
\node[draw,circle,minimum size=\dia,inner sep=0](b24) at ({90 + (\Ang * 2) + \anga}:\radb) {\footnotesize $b_{24}$};
\node[draw,circle,minimum size=\dia,inner sep=0](b7) at ({90 + (\Ang * 2) }:\radb) {\footnotesize $b_{7}$};
\node[draw,circle,minimum size=\dia,inner sep=0](b25) at ({90 + (\Ang * 2) - \anga}:\radb) {\footnotesize $b_{25}$};
\node[draw,circle,minimum size=\dia,inner sep=0](a24) at ({90 + (\Ang * 2) + \anga}:\rada) {\footnotesize $a_{24}$};
\node[draw,circle,minimum size=\dia,inner sep=0](a25) at ({90 + (\Ang * 2) - \anga}:\rada) {\footnotesize $a_{25}$};
\draw[color=blue] (a1) -- (a2) -- (a3) -- (a4) -- (a5) -- (a6) -- (a7) -- (a8) -- (a9) -- (a10) -- (a11) -- (a12) -- (a13) -- (a14) -- (a15) -- (a16) -- (a17) -- (a18) -- (a19) -- (a20) -- (a21) -- (a22) -- (a23) -- (a24) -- (a25) -- (a1);
\draw[color=gren] (b1) -- (b3) -- (b20) -- (b2) -- (b19) -- (b1);
\draw[color=gren] (b25) -- (b7) -- (b24) -- (b6) -- (b8) -- (b25);
\draw[color=gren] (b4) -- (b11) -- (b13) -- (b5) -- (b12) -- (b4);
\draw[color=gren] (b23) -- (b21) -- (b14) -- (b22) -- (b15) -- (b23);
\draw[color=gren] (b18) -- (b10) -- (b17) -- (b9) -- (b16) -- (b18);
\draw[color=red] (a1) -- (b1);
\draw[color=red] (a2) -- (b2);
\draw[color=red] (a3) -- (b3);
\draw[color=red] (a4) -- (b4);
\draw[color=red] (a5) -- (b5);
\draw[color=red] (a6) -- (b6);
\draw[color=red] (a7) -- (b7);
\draw[color=red] (a8) -- (b8);
\draw[color=red] (a9) -- (b9);
\draw[color=red] (a10) -- (b10);
\draw[color=red] (a11) -- (b11);
\draw[color=red] (a12) -- (b12);
\draw[color=red] (a13) -- (b13);
\draw[color=red] (a14) -- (b14);
\draw[color=red] (a15) -- (b15);
\draw[color=red] (a16) -- (b16);
\draw[color=red] (a17) -- (b17);
\draw[color=red] (a18) -- (b18);
\draw[color=red] (a19) -- (b19);
\draw[color=red] (a20) -- (b20);
\draw[color=red] (a21) -- (b21);
\draw[color=red] (a22) -- (b22);
\draw[color=red] (a23) -- (b23);
\draw[color=red] (a24) -- (b24);
\draw[color=red] (a25) -- (b25);
\end{tikzpicture}
\caption{The Watkins graph, shown in its traditional drawing.}
\label{f:WatkinsUsual}
\end{figure}
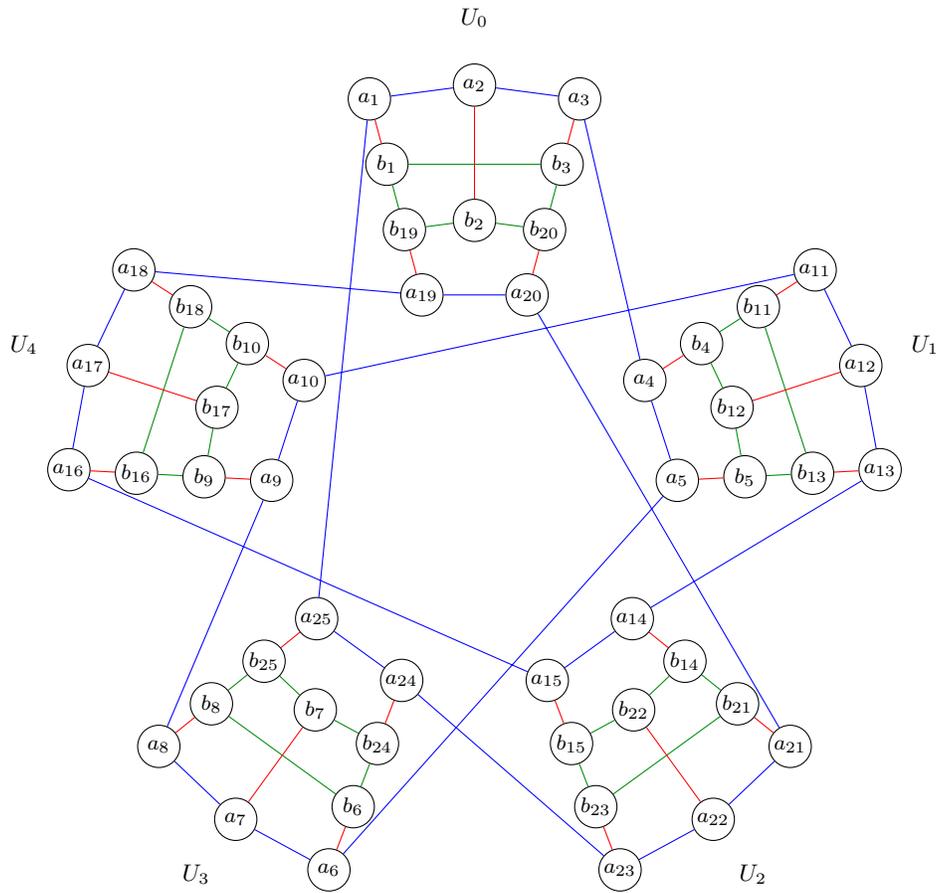

\begin{figure}
\centering
\begin{tikzpicture}[scale=0.65]
\def \lev {3.0}
\def \eps {.5}
\def \wid {20}
\def \ext {1}
\def \dia {20}
\def \thk {1.4}
\def \thn {1.1}
\node[draw=none,fill=none](L0) at (0,{0*\lev+2*\eps}) {$V_0$};
\node[draw,red,line width=\thk,circle,minimum size=\dia,inner sep=0](a1) at ({1*\wid/2},{0*\lev+2*\eps}) {$a_{1}$};
\node[draw=none,fill=none](L1) at (0,{1*\lev+\eps}) {$V_1$};
\node[draw,blue,line width=\thk,circle,minimum size=\dia,inner sep=0](a2) at ({1*\wid/4},{1*\lev+\eps}) {$a_{2}$};
\node[draw,rectangle,minimum size=\dia,inner sep=0](b1) at ({2*\wid/4},{1*\lev+\eps}) {$b_{1}$};
\node[draw,gren,line width=\thk,circle,minimum size=\dia,inner sep=0](a25) at ({3*\wid/4},{1*\lev+\eps}) {$a_{25}$};
\node[draw=none,fill=none](L2) at (0,{2*\lev}) {$V_2$};
\node[draw,blue,line width=\thk,circle,minimum size=\dia,inner sep=0](a3) at ({1*\wid/7},{2*\lev}) {$a_{3}$};
\node[draw,rectangle,minimum size=\dia,inner sep=0](b3) at ({2*\wid/7},{2*\lev}) {$b_{3}$};
\node[draw,rectangle,minimum size=\dia,inner sep=0](b2) at ({3*\wid/7},{2*\lev}) {$b_{2}$};
\node[draw,rectangle,minimum size=\dia,inner sep=0](b19) at ({4*\wid/7},{2*\lev}) {$b_{19}$};
\node[draw,rectangle,minimum size=\dia,inner sep=0](b25) at ({5*\wid/7},{2*\lev}) {$b_{25}$};
\node[draw,gren,line width=\thk,circle,minimum size=\dia,inner sep=0](a24) at ({6*\wid/7},{2*\lev}) {$a_{24}$};
\node[draw=none,fill=none](L3) at (0,{3*\lev}) {$V_3$};
\node[draw,blue,line width=\thk,circle,minimum size=\dia,inner sep=0](a4) at ({1*\wid/8},{3*\lev}) {$a_{4}$};
\node[draw,rectangle,minimum size=\dia,inner sep=0](b20) at ({2*\wid/8},{3*\lev}) {$b_{20}$};
\node[draw,rectangle,minimum size=\dia,inner sep=0](b7) at ({3*\wid/8},{3*\lev}) {$b_{7}$};
\node[draw,rectangle,minimum size=\dia,inner sep=0](b24) at ({4*\wid/8},{3*\lev}) {$b_{24}$};
\node[draw,rectangle,minimum size=\dia,inner sep=0](b8) at ({5*\wid/8},{3*\lev}) {$b_{8}$};
\node[draw,circle,minimum size=\dia,inner sep=0](a19) at ({6*\wid/8},{3*\lev}) {$a_{19}$};
\node[draw,gren,line width=\thk,circle,minimum size=\dia,inner sep=0](a23) at ({7*\wid/8},{3*\lev}) {$a_{23}$};
\node[draw=none,fill=none](L4) at (0,{4*\lev}) {$V_4$};
\node[draw,blue,line width=\thk,rectangle,minimum size=\dia,inner sep=0](b4) at ({1*\wid/10},{4*\lev}) {$b_{4}$};
\node[draw,blue,line width=\thk,circle,minimum size=\dia,inner sep=0](a5) at ({2*\wid/10},{4*\lev}) {$a_{5}$};
\node[draw,text=blue,line width=\thn,circle,minimum size=\dia,inner sep=0](a7) at ({3*\wid/10},{4*\lev}) {$a_{7}$};
\node[draw,circle,minimum size=\dia,inner sep=0](a8) at ({4*\wid/10},{4*\lev}) {$a_{8}$};
\node[draw,circle,minimum size=\dia,inner sep=0](a18) at ({5*\wid/10},{4*\lev}) {$a_{18}$};
\node[draw,text=blue,line width=\thn,rectangle,minimum size=\dia,inner sep=0](b6) at ({6*\wid/10},{4*\lev}) {$b_{6}$};
\node[draw,text=gren,line width=\thn,circle,minimum size=\dia,inner sep=0](a20) at ({7*\wid/10},{4*\lev}) {$a_{20}$};
\node[draw,gren,line width=\thk,circle,minimum size=\dia,inner sep=0](a22) at ({8*\wid/10},{4*\lev}) {$a_{22}$};
\node[draw,gren,line width=\thk,rectangle,minimum size=\dia,inner sep=0](b23) at ({9*\wid/10},{4*\lev}) {$b_{23}$};
\node[draw=none,fill=none,text=blue] at ({3*\wid/10},{4*\lev+\ext}) {$1$};
\node[draw=none,fill=none,text=blue] at ({6*\wid/10},{4*\lev+\ext}) {$1$};
\node[draw=none,fill=none,text=gren] at ({7*\wid/10},{4*\lev+\ext}) {$1$};
\node[draw=none,fill=none](L5) at (0,{5*\lev}) {$V_5$};
\node[draw,blue,line width=\thk,rectangle,minimum size=\dia,inner sep=0](b11) at ({1*\wid/12},{5*\lev}) {$b_{11}$};
\node[draw,blue,line width=\thk,rectangle,minimum size=\dia,inner sep=0](b5) at ({2*\wid/12},{5*\lev}) {$b_{5}$};
\node[draw,blue,line width=\thk,rectangle,minimum size=\dia,inner sep=0](b12) at ({3*\wid/12},{5*\lev}) {$b_{12}$};
\node[draw,text=blue,line width=\thn,circle,minimum size=\dia,inner sep=0](a6) at ({4*\wid/12},{5*\lev}) {$a_{6}$};
\node[draw,circle,minimum size=\dia,inner sep=0](a9) at ({5*\wid/12},{5*\lev}) {$a_{9}$};
\node[draw,rectangle,minimum size=\dia,inner sep=0](b18) at ({6*\wid/12},{5*\lev}) {$b_{18}$};
\node[draw,circle,minimum size=\dia,inner sep=0](a17) at ({7*\wid/12},{5*\lev}) {$a_{17}$};
\node[draw,text=gren,line width=\thn,circle,minimum size=\dia,inner sep=0](a21) at ({8*\wid/12},{5*\lev}) {$a_{21}$};
\node[draw,gren,line width=\thk,rectangle,minimum size=\dia,inner sep=0](b21) at ({9*\wid/12},{5*\lev}) {$b_{21}$};
\node[draw,gren,line width=\thk,rectangle,minimum size=\dia,inner sep=0](b22) at ({10*\wid/12},{5*\lev}) {$b_{22}$};
\node[draw,gren,line width=\thk,rectangle,minimum size=\dia,inner sep=0](b15) at ({11*\wid/12},{5*\lev}) {$b_{15}$};
\node[draw=none,fill=none,text=blue] at ({4*\wid/12},{5*\lev+\ext}) {$1$};
\node[draw=none,fill=none,text=gren] at ({8*\wid/12},{5*\lev+\ext}) {$1$};
\node[draw=none,fill=none](L6) at (0,{6*\lev}) {$V_6$};
\node[draw,blue,line width=\thk,rectangle,minimum size=\dia,inner sep=0](b13) at ({1*\wid/12},{6*\lev}) {$b_{13}$};
\node[draw,blue,line width=\thk,circle,minimum size=\dia,inner sep=0](a12) at ({2*\wid/12},{6*\lev}) {$a_{12}$};
\node[draw,text=oran,line width=\thn,circle,minimum size=\dia,inner sep=0](a11) at ({3*\wid/12},{6*\lev}) {$a_{11}$};
\node[draw,text=oran,line width=\thn,circle,minimum size=\dia,inner sep=0](a10) at ({4*\wid/12},{6*\lev}) {$a_{10}$};
\node[draw,text=oran,line width=\thn,rectangle,minimum size=\dia,inner sep=0](b10) at ({5*\wid/12},{6*\lev}) {$b_{10}$};
\node[draw,text=oran,line width=\thn,rectangle,minimum size=\dia,inner sep=0](b17) at ({6*\wid/12},{6*\lev}) {$b_{17}$};
\node[draw,text=oran,line width=\thn,rectangle,minimum size=\dia,inner sep=0](b9) at ({7*\wid/12},{6*\lev}) {$b_{9}$};
\node[draw,text=oran,line width=\thn,rectangle,minimum size=\dia,inner sep=0](b16) at ({8*\wid/12},{6*\lev}) {$b_{16}$};
\node[draw,text=oran,line width=\thn,circle,minimum size=\dia,inner sep=0](a16) at ({9*\wid/12},{6*\lev}) {$a_{16}$};
\node[draw,gren,line width=\thk,circle,minimum size=\dia,inner sep=0](a15) at ({10*\wid/12},{6*\lev}) {$a_{15}$};
\node[draw,gren,line width=\thk,rectangle,minimum size=\dia,inner sep=0](b14) at ({11*\wid/12},{6*\lev}) {$b_{14}$};
\node[draw=none,fill=none,text=oran] at ({3*\wid/12},{6*\lev+\ext}) {$1$};
\node[draw=none,fill=none,text=oran] at ({4*\wid/12},{6*\lev+\ext}) {$1$};
\node[draw=none,fill=none,text=oran] at ({5*\wid/12},{6*\lev+\ext}) {$1$};
\node[draw=none,fill=none,text=oran] at ({6*\wid/12},{6*\lev+\ext}) {$1$};
\node[draw=none,fill=none,text=oran] at ({7*\wid/12},{6*\lev+\ext}) {$1$};
\node[draw=none,fill=none,text=oran] at ({8*\wid/12},{6*\lev+\ext}) {$1$};
\node[draw=none,fill=none,text=oran] at ({9*\wid/12},{6*\lev+\ext}) {$1$};
\node[draw=none,fill=none](L7) at (0,{7*\lev}) {$V_7$};
\node[draw,blue,line width=\thk,circle,minimum size=\dia,inner sep=0](a13) at ({1*\wid/3},{7*\lev}) {$a_{13}$};
\node[draw,gren,line width=\thk,circle,minimum size=\dia,inner sep=0](a14) at ({2*\wid/3},{7*\lev}) {$a_{14}$};
\node[draw=none,fill=none] at ({1*\wid/3},{7*\lev+\ext}) {$85$};
\node[draw=none,fill=none] at ({2*\wid/3},{7*\lev+\ext}) {$85$};
\draw (a1) -- (a2) -- (a3) -- (a4) -- (a5) -- (a6) -- (a7) -- (a8) -- (a9) -- (a10) -- (a11) -- (a12) -- (a13) -- (a14) -- (a15) -- (a16) -- (a17) -- (a18) -- (a19) -- (a20) -- (a21) -- (a22) -- (a23) -- (a24) -- (a25) -- (a1);
\draw (b1) -- (b3) -- (b20) -- (b2) -- (b19) -- (b1);
\draw (b25) -- (b7) -- (b24) -- (b6) -- (b8) -- (b25);
\draw (b4) -- (b11) -- (b13) -- (b5) -- (b12) -- (b4);
\draw (b23) -- (b21) -- (b14) -- (b22) -- (b15) -- (b23);
\draw (b18) -- (b10) -- (b17) -- (b9) -- (b16) -- (b18);
\draw (a1) -- (b1);
\draw (a2) -- (b2);
\draw (a3) -- (b3);
\draw (a4) -- (b4);
\draw (a5) -- (b5);
\draw (a6) -- (b6);
\draw (a7) -- (b7);
\draw (a8) -- (b8);
\draw (a9) -- (b9);
\draw (a10) -- (b10);
\draw (a11) -- (b11);
\draw (a12) -- (b12);
\draw (a13) -- (b13);
\draw (a14) -- (b14);
\draw (a15) -- (b15);
\draw (a16) -- (b16);
\draw (a17) -- (b17);
\draw (a18) -- (b18);
\draw (a19) -- (b19);
\draw (a20) -- (b20);
\draw (a21) -- (b21);
\draw (a22) -- (b22);
\draw (a23) -- (b23);
\draw (a24) -- (b24);
\draw (a25) -- (b25);
\end{tikzpicture}
\caption{The Watkins graph, organized by distance from target $a_1$, along with an $a_1$-unsolvable configuration $C$ of size $182$.}
\label{f:WatkinsDista1}
\end{figure}
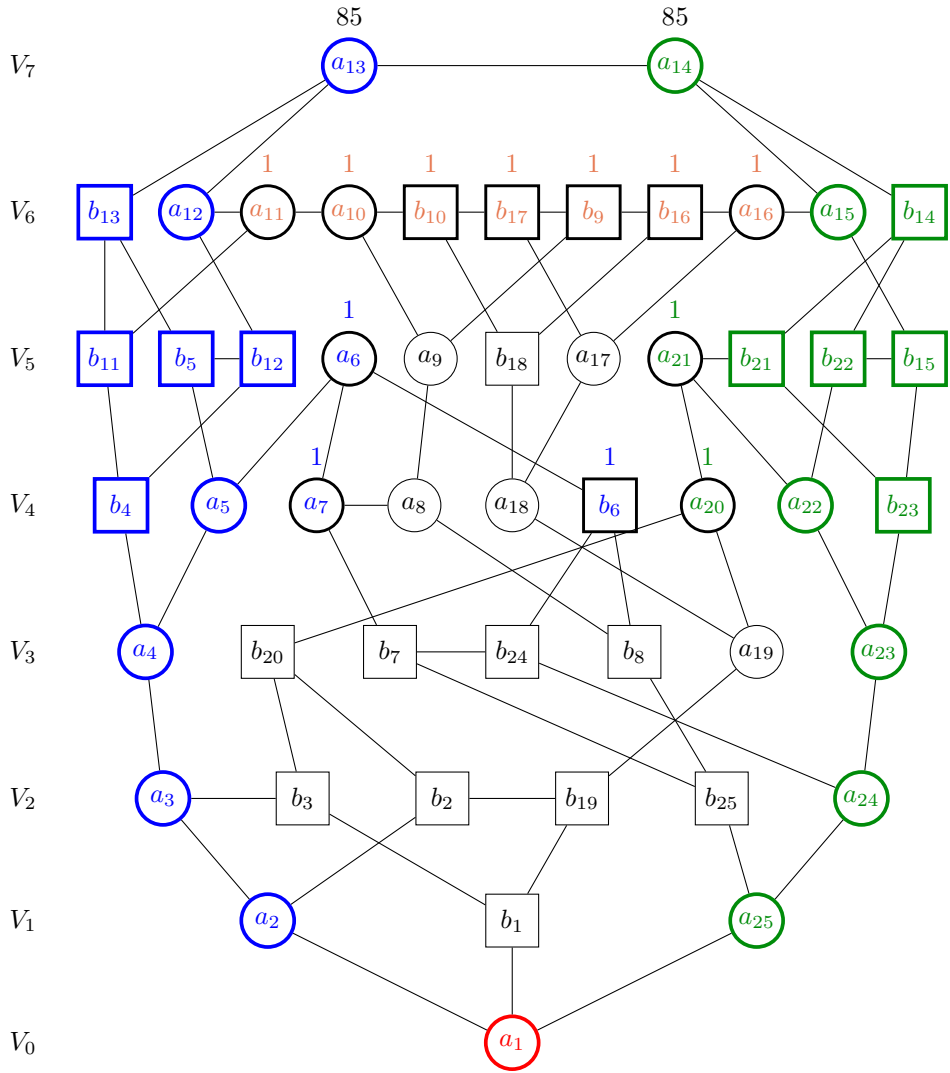

We first observe an isomorphism between the traditional drawing of $W$ shown in Figure \ref{f:WatkinsUsual} and the drawing of $W$ in Figure \ref{f:WatkinsDista1} that partitions its vertices by distance from $a_1$ (those of distance $i$ are in $V_i$).
Indeed, one can check that, in both drawings, $W$ contains the $25$-cycle $a_1, a_2, \ldots, a_{25}, a_1$ (shown with blue edges in Figure \ref{f:WatkinsUsual} and circular vertex shapes in Figure \ref{f:WatkinsDista1}), five vertex-disjoint $5$-cycles on the vertices $b_1, \ldots, b_{25}$ (shown with green edges in Figure \ref{f:WatkinsUsual} and square vertex shapes in Figure \ref{f:WatkinsDista1}), and the perfect matching using the $25$ pairs $\{a_i,b_i\}$ (shown with red edges in Figure \ref{f:WatkinsUsual}).
The point of showing Figure \ref{f:WatkinsUsual} is to recognize $W$ and its symmetries; the point of showing Figure~\ref{f:WatkinsDista1} is to make lower bound arguments about $\pi(W)$.

Next we prove the lower bound by displaying an $a_1$-unsolvable configuration of size $182$; such a configuration $C$ is shown in Figure~\ref{f:WatkinsDista1}, with each vertex of $V_7$ having $85$ pebbles and 12 other vertices (7 orange, 3 blue, and 2 green) having one pebble each.
The argument we give will make use of the vertex colors in Figure \ref{f:WatkinsDista1}.

For two vertices $u$ and $v$, we call a $uv$-path of minimal length a $uv$-{\it geodesic}.
Define $V_B$ (resp. $V_G$) to be the set of vertices in the union of all $a_{13}a_1$-geodesics (resp. $a_{14}a_1$-geodesics); these are shown in Figure \ref{f:WatkinsDista1} with blue (resp. green) shapes, except for the target vertex $a_1$ in red.
Note that the path $P_8$ is a retract of the subgraph of $W$ induced by $V_B$ (resp. $V_G$), and so the graph $W_{BG}$ induced by $V_B\cup V_G$ has the cycle $C_{15}$ as a retract.
Hence $\pi(W_{BG})\ge \pi(C_{15})=171$.
In particular, by Lemma \ref{l:retract}, the configuration $C_{BG}$ --- the restriction of the configuration $C$ shown in Figure \ref{f:WatkinsDista1} to $V_B\cup V_G$ --- is $a_1$-unsolvable in $W_{BG}$.

Suppose that $C_{BG}$ has a minimal $a_1$-solution $\sigma$ in $W$.
If $\sigma$ uses the edge $e_2 = a_{2}a_1$ then it does not use either of the edges $b_1a_1$ or $a_{25}a_1$; hence $\sigma$ is an $a_1$-solution in the subgraph $W_2 = W-e_2$.
However the distance from $a_{14}$ to $a_2$ is $7$, and so the $a_1$-weight of $C_{BG}$ equals $85(3)/2^8 = 255/256<1$, a contradiction.
Similarly, if $\sigma$ uses the edge $e_1=b_1a_1$ then we obtain the same $a_1$-weight calculation, while if $\sigma$ uses the edge $e_{25}=a_{25}a_1$ then we obtain a smaller $a_1$-weight calculation; in both cases we have a contradiction.
Hence $C_{BG}$ is not $a_1$-solvable in $W$.
(This already invalidates the claim of \cite{pebblingwatkinssnark} that $\pi(W)=166$.)
Therefore, if $C$ is $a_1$-solvable in $W$, then any solution must use some of the pebbles in the configuration $S = C - C_{BG}$ of singletons of $C$.

Now let $\sigma$ be any minimal $a_1$-solution and let $S_\sigma$ denote the set of singletons of $S$ used by $\sigma$.
We may rearrange the order of pebbling moves so that those moves involving $S_\sigma$ are performed as early as possible.
Now let $C_\sigma$ denote the configuration that results from halting future pebbling moves of $\sigma$ once the moves involving $S_\sigma$ have been performed.
Additionally, denote by $C'_\sigma$ the configuration $C_\sigma-(S-S_\sigma)$; that is, from $C_\sigma$, throw away the singletons of $S$ that were not used by $\sigma$.
Obviously $C'_\sigma$ is $a_1$-solvable.
We argue that this is impossible by showing that the $a_1$-weight of $C'_\sigma$ is less than 1.

A $uv$-{\it slide} is a path $ua_1\cdots a_kv$ in which $u$ is big and each $a_i$ is nonempty, the use of which allows for a pebble to move from $u$ to $v$ along the slide.
A {\it potential slide} ignores the requirement that $u$ is big; it is ``potential'' because it becomes a slide if $u$ becomes big.
Notice that every singleton vertex is part of a potential slide.
Because the $a_1$-weight of $C_{BG}$ in $W$ equals $1-2^{-8}$, $C'_\sigma$ can only be $a_1$-solvable if we are able to increase the $a_1$-weight of $C_{BG}$ by at least $2^{-8}$ by using the slides of $S$.
In other words, it must be that the $a_1$-weight of a pebble in $C'_\sigma$ must be greater than the sum of the $a_1$-weight of the pebbles of $C_{BG}$ that were used to create it.
We show that this cannot happen.

Partition $S=S_6\cup S_B\cup S_G$, where $S_6=S\cap V_6$ (the orange pebbles in Figure \ref{f:WatkinsDista1}), $S_B=S\cap V_B$ (the blue pebbles in Figure \ref{f:WatkinsDista1}), and $S_G=S\cap V_G$ (the green pebbles in Figure \ref{f:WatkinsDista1}).
Using a slide in $S_6$ requires two pebbles from $V_5\cup V_6$ and places one pebble in either $V_5\cup V_6$. 
Such a pebble has the same or lesser $a_1$-weight than the total $a_1$-weight of the original two pebbles.
Using a slide in $S_B$ requires two pebbles from $V_3\cup V_4$ and places one pebble in $V_3\cup V_4$; this also cannot increase the $a_1$-weight.
Using a slide in $S_G$, however, can increase the $a_1$-weight only if the two pebbles originate from $b_{21}$ and the resulting pebble lands on $b_{20}$.
Still, the resulting pebble must be used in $\sigma$ to achieve this $a_1$-weight increase (of $1/2^3-2/2^5 = 1/2^4$), although to do so requires at least $2^6$ pebbles from $V_7$, resulting in a $a_1$-weight loss of $2^7/2^7-1/2^3 > 1/2^4$.
Hence it is impossible to increase the $a_1$-weight of $C'_\sigma$ via pebbling steps, and so $C'_\sigma$, and thus $C$, is $a_1$-unsolvable.
\hfill $\Box$

\section{Final remarks}
\label{s:final}

The previous work of \cite{pebblingflowersnark} on Flower snarks and of \cite{pebblingwatkinssnark} on the Watkins snark 
are at the foundation of our findings, since after reading their published papers, we have realised we were able to contribute to the subject.
We hope the developed techniques and results will motivate more researchers to achieve better bounds on the pebbling numbers of snarks.
Table 1 shows the state of art of the pebbling numbers of several well known snarks, using the basic bounds mentioned in the introduction, as well as Lemma~\ref{l:SNL}, and Theorems~\ref{t:Petersen},~\ref{t:Flower},~\ref{t:Blanusa} and~\ref{t:Watkins}.

\begin{table}
\resizebox{\textwidth}{!}{
\begin{tabular}{|l|c|c|c|}
\hline
Snark & $n(G)$ & $D(G)$ & $\pi(G)$ \\
\hline
Petersen & 10 & 2 & 10 \\ 
\hline
Flower $J_3$ & 12 & 3 & $12 \le \pi(J_3)\le 13$ \\
\hline
Flower $J_5$ & 20 & 4 & $23\le \pi(J_5)\le 30$ \\
\hline
Flower $J_7$ & 28 & 5 & $41\le \pi(J_7)\le 61$ \\
\hline
Flower $J_m$ ($m = 2k+1 \ge 7$) & $4m$ & $k+2$ & $2^{k+2}+9$
$\le \pi(J_m)\le \lfloor 2^{k+2}9/5 + 2k - 18/5 \rfloor +1$ \\
\hline
Blanu\v{s}a (1 and 2) & 18 & 4 & $ 23 \le \pi(G) \leq 34$ \\
\hline
Loupekine (1 and 2) & 22 & 4 & $ 24 \le \pi(G) \leq 271$ \\
\hline
Double-Star & 30 & 4 & $32 \le \pi(G) \leq 391$ \\
\hline
Szekeres & 50 & 7 & $128 \le \pi(G) \leq 5462$ \\
\hline
Watkins & 50 & 7 & $183 \le \pi(G) \leq 5462$ \\
\hline
 \end{tabular}}
\caption{Bounds on the pebbling numbers of several well known snarks.}
\label{snarktable}
\end{table}

\section*{Acknowledgements}
This research was partially supported by CAPES  Finance Code 001, FAPERJ (grant no. E-26/202.793/2017, E-26/010.002239/201, E-26/010.002674/2019, E-26/201.360/2021), and CNPq (grant no. 407635/2018-1, 305356/2021-6, 313797/2020-0).
The third author is grateful to the Federal University of Rio de Janeiro for partial support of his visit, and to his coauthors for their tremendous hospitality during his visit.

\bibliographystyle{acm}
\bibliography{refs}

\section*{Appendix}

\subsection*{Missing cases for the proof of Theorem~\ref{t:Flower}}
Regarding the upper bounds of Theorem~\ref{t:Flower}, for the two other possible roots $v_0$ and $x_0$, we have obtained below strategies that give values not larger than the ones given by strategies presented for root $z_0$.

For $J_3$, we define three $v_0$-strategies $\bT_0$, $\bT_1$, and $\bT_{-1}$ by 
\vspace{-0.04in}
\begin{quote}
    $\bullet$ 
    $\bT_0(z_0,x_0,y_0,x_1,y_1,x_{-1},y_{-1},z_1,z_{-1}) = (8,4,4,2,2,2,2,1,1)$, \\
    $\bullet$
    $\bT_1(v_1,z_1,x_1,y_1,x_0,x_{-1},y_{-1}) = (8,4,2,2,1,1,1)$ and\\ 
    $\bullet$
    $\bT_{-1}(v_{-1},z_{-1},x_{-1},y_{-1},y_0,x_1,y_1) = (8,4,2,2,1,1,1)$,
\end{quote}
giving rise to the inequality $|C| \le \frac{1}{5}(\bT_0+\bT_1+\bT_{-1}) \le 64/5$ whenever $C$ is $v_0$-unsolvable.
Hence $\pi(J_3,v_0) \le 13$.

And we define three $x_0$-strategies $\bT_0$, $\bT_1$, and $\bT_{-1}$ by 
\vspace{-0.04in}
\begin{quote}
    $\bullet$ 
    $\bT_0(z_0,v_0,y_0,v_{-1},v_1,y_{-1},y_1) = (8,4,4,2,2,1,1)$, \\
    $\bullet$
    $\bT_1(x_1,z_1,y_{-1},v_1,z_{-1},y_0,v_{-1},v_0) = (8,4,4,2,1,1,1,1)$ and\\ 
    $\bullet$
    $\bT_{-1}(x_{-1},z_{-1}, y_1, v_{-1,}z_1,v_1) = (8,4,4,2,1,1)$,
\end{quote}
giving rise again to the inequality $|C| \le \frac{1}{5}(\bT_0+\bT_1+\bT_{-1}) \le 64/5$ whenever $C$ is $x_0$-unsolvable.
Hence $\pi(J_3,x_0) \le 13$.
We may conclude that $\pi(J_3) \le 13$.

For $J_5$, we define three $v_0$-strategies $\bT_0$, $\bT_1$, and $\bT_{-1}$ by 
\vspace{-0.04in}
\begin{quote}
    $\bullet$ 
    $\bT_0(z_0,x_0,y_0,x_1,y_1,x_{-1},y_{-1},x_2,y_2,x_{-2},y_{-2},z_2,z_{-2}) = (16,8,8,4,4,4,4,2,2,1,1,1,1)$, \\
    $\bullet$
    $\bT_1(v_1,v_2, z_1, z_2, x_2,y_2,x_1,y_{-2},x_{-2},y_1) =  
    (16,8,8,4,2,2,1,1,1,1)$ and\\ 
    $\bullet$
    $\bT_{-1}(v_{-1},v_{-2}, z_{-1}, z_{-2},x_{-2},y_{-2},x_{-1}, y_2,   x_2, y_{-1}) = (16,8,8,4,2,2,1,1,1,1)$,
\end{quote}
giving rise to the inequality $|C| \le \frac{1}{5}(\bT_0+\bT_1+\bT_{-1}) \le 146/5$ whenever $C$ is $v_0$-unsolvable.
Hence $\pi(J_5,v_0) \le 30$.

And we define three $x_0$-strategies $\bT_0$, $\bT_1$, and $\bT_{-1}$ by 
\vspace{-0.04in}
\begin{quote}
    $\bullet$ 
    $\bT_0(z_0,v_0,y_0,v_{-1},v_1,y_{-1},y_1,v_2,v_{-2},y_2,y_{-2},z_2) = (16,8,8,4,4,4,4,2,2,1,1,1)$, \\
    $\bullet$
    $\bT_1(x_1,z_1,x_2,y_{-2},z_2,  y_{-1},z_{-2},y_2,v_2, v_1,v_{-2}) = (16,8,8,4,4,1,1,1,2,1,1)$ and\\ 
    $\bullet$
    $\bT_{-1}(x_{-1},z_{-1},x_{-2},z_{-2}, y_2, v_{-2},y_1,v_2,v_{-1}) = (16,8,8,4,4,2,1,1,1)$,
\end{quote}
giving rise again to the inequality $|C| \le \frac{1}{5}(\bT_0+\bT_1+\bT_{-1}) \le 146/5$ whenever $C$ is $x_0$-unsolvable.
Hence $\pi(J_5,x_0) \le 30$.
We may conclude that $\pi(J_5) \le 30$.

For $J_7$, we define three $v_0$-strategies $\bT_0$, $\bT_1$, and $\bT_{-1}$ by 
\vspace{-0.04in}
\begin{quote}
    $\bullet$ 
    $\bT_0(z_0,x_0,y_0,x_1,y_1,x_{-1},y_{-1},x_2,z_1,z_{-1},x_{-2},y_2,y_{-2},x_3,x_{-3}, y_3,y_{-3},z_3,z_{-3})$ \\
    $= (32,16,16,8,8,8,8,4,4,4,4,4,4,2,2,2,2,1,1)$, \\
    $\bullet$
    $\bT_1(v_1,v_2, z_1, v_3,z_2,z_3, x_2,y_2,x_3,y_3,x_{-3},y_{-3}) =  
    (32,16,1,8,5,4,1,1,2,2,1,1)$ and\\ 
    $\bullet$
    $\bT_{-1}(v_{-1},v_{-2}, z_{-1}, v_{-3},z_{-2},z_{-3},x_{-2},y_{-2},x_{-3},y_{-3}, y_3,   x_3) = (32,16,1,8,5,4,1,1,2,2,1,1)$,
\end{quote}
giving rise to the inequality $|C| \le \frac{1}{5}(\bT_0+\bT_1+\bT_{-1}) \le 278/5 <56$ whenever $C$ is $v_0$-unsolvable.
Hence $\pi(J_7,v_0) \le 56$.

And we define three $x_0$-strategies $\bT_0$, $\bT_1$, and $\bT_{-1}$ by 
\vspace{-0.04in}
\begin{quote}
    $\bullet$ 
    $\bT_0(z_0,v_0,y_0,v_1,v_{-1},y_1,y_{-1},v_2,v_{-2},y_2,y_{-2},v_3,v_{-3},y_3,y_{-3},z_3)$ \\
    $= (32,16,16,8,8,8,8,4,4,4,4,2,2,1,1,1)$, \\
    $\bullet$
    $\bT_1(x_1,z_1,x_2,z_2,x_3,v_2,   y_{-3},z_3,z_{-3},   y_3, v_3,v_{-3}) = (32,5,16,8,8,1,4,4,1,1,2,1)$ and\\ 
    $\bullet$
    $\bT_{-1}(x_{-1},z_{-1}, x_{-2}, z_{-2}, x_{-3}, y_{-2}, v_{-2}, z_{-3},   y_3,   v_{-3},y_2,v_3) = (32,5,16,8,8,1,1,4,4,2,1,1)$,
\end{quote}
giving rise to the inequality $|C| \le \frac{1}{5}(\bT_0+\bT_1+\bT_{-1}) \le 284/5 <57$ whenever $C$ is $x_0$-unsolvable.
Hence $\pi(J_7,x_0) \le 57$.
We may conclude that $\pi(J_7) \le 61$.

For $J_9$, we define three $x_0$-strategies $\bT_0$, $\bT_1$, and $\bT_{-1}$ by 
\vspace{-0.04in}
\begin{quote}
    $\bullet$ 
    $\bT_0(z_0,v_0,y_0,v_1,v_{-1},y_1,y_{-1},v_2,v_{-2},y_2,y_{-2}, v_3,v_{-3},y_3,y_{-3}, v_4,v_{-4},y_4,y_{-4},z_4)$ \\
    $= (64,32,32,16,16,16,16,8,8,8,8,4,4,4,4,2,2,2,1,1)$, \\
    $\bullet$
    $\bT_1(x_1,z_1,x_2,z_2,x_3,  z_3,x_4, v_3,  y_{-4},z_4,z_{-4},   y_4, v_4,v_{-4}) = (64,5,32,16,16, 5, 8,1,2,4,1,1,2,1)$ and\\ 
    $\bullet$
    $\bT_{-1}(x_{-1},z_{-1},x_{-2},z_{-2}, x_{-3},v_{-2}, z_{-3}, x_{-4},y_{-3}, v_{-3}, z_{-4},y_4,v_{-4},y_3,v_4)$ \\
    $= (64,5,32,16,16,8,5,8,1,1,4,2,2,1,1)$.
\end{quote}

For $m\ge 7$ (i.e. $k\ge 3$), using the same pattern we have defined above for the three $v_0$-strategies for $J_7$, we define three corresponding $v_0$-strategies
giving rise to the inequality 
\begin{align*}
    5|C|
    &\le \bT_0+\bT_1+\bT_{-1}\\
  & = 3(2^{k+2}) +   4(2^{k+1})   +6(2^{k} + \cdots + 2^3)  + 5(2k+8) \\
  &\le 3(2^{k+2}) + 6(2^3 + \cdots + 2^{k+1})  + 5(2k+6),
\end{align*}
whenever $C$ is $v_0$-unsolvable.

For $m\ge 9$ (i.e. $k\ge 4$), using the same pattern we have defined above for the three $x_0$-strategies for $J_9$, we define three corresponding $x_0$-strategies
giving rise to the inequality 
\begin{align*}
    5|C|
    &\le \bT_0+\bT_1+\bT_{-1}\\
  & = 3(2^{k+2}) +   4(2^{k+1})   + 8(2^{k})  +6(2^{k-1} +  \cdots + 2^3) + 5(2k+6) \\
  &\le 3(2^{k+2}) + 6(2^3 + \cdots + 2^{k+1})  + 5(2k+6),
\end{align*}
whenever $C$ is $x_0$-unsolvable.

We may conclude that for $m\ge 3$ (i.e. $k\ge 1$), 
$\pi(J_m)\le \lfloor 2^{k+2}9/5 + 2k - 18/5 \rfloor +1$.

\subsection*{Missing cases for the proof of Theorem~\ref{t:Blanusa}}

Regarding the upper bounds of Theorem~\ref{t:Blanusa}, for the five other possible roots $z_1$, $x_1$, $z_2$, $x_2$, and $z_3$, we have obtained below strategies that give values not larger than the ones given by strategies presented for root $x_3$.

For $B_2$, we define three $x_1$-strategies $\bT_1$, $\bT_2$, and $\bT_3$ by 
\vspace{-0.04in}
\begin{quote}
    $\bullet$ 
    $\bT_1(x_3,z_3,x'_3,x'_1,x'_5,z'_1) = (32,16,16,8,8,4)$, \\
    $\bullet$
    $\bT_2(z_1,z_2,z_5,x_2,x_5,z'_2,z'_5,x'_2,x'_5,z'_1) = (16,8,8,4,4,4,4,2,2,2)$ and\\ 
    $\bullet$
    $\bT_3(x_4,z_4, x'_4,x'_1,x'_2) = (32,16,8,8,4$)
\end{quote}
giving rise to the inequality $|C| \le \frac{1}{4}(\bT_1+\bT_2+\bT_3) \le 120/4$ whenever $C$ is $x_1$-unsolvable.
Hence $\pi(B_2,x_1)\le 31$.

For $B_2$, we define three $x_2$-strategies $\bT_1$, $\bT_2$, and $\bT_3$ by 
\vspace{-0.04in}
\begin{quote}
    $\bullet$ 
    $\bT_1(x_4,z_4,x'_4,x'_1,x'_3) = (16,8,4,2,1)$, \\
    $\bullet$
    $\bT_2(z_2,z_1,z'_2,x_1,x'_2,z'_1,x'_1,x'_5,x'_3) = (16,8,8,4,4,4,2,2,1)$ and\\ 
    $\bullet$
    $\bT_3(x_5,x_3, z_5,z_3,z'_5,x'_3,x'_5) = (16,8,4,4,2,2$)
\end{quote}
giving rise to the inequality $|C| \le \frac{1}{4}(\bT_1+\bT_2+\bT_3) \le 124/4$ whenever $C$ is $x_2$-unsolvable.
Hence $\pi(B_2,x_2)\le 32$.

For $B_2$, we define three $z_1$-strategies $\bT_1$, $\bT_2$, and $\bT_3$ by 
\vspace{-0.04in}
\begin{quote}
    $\bullet$ 
    $\bT_1(x_1,x_3,x_4,z_3,z_4,x'_3,x'_4,x'_1) = (16,8,8,4,4,2,2,1)$, \\
    $\bullet$
    $\bT_2(z_2,x_2,z'_2,x'_1,z'_1,x'_4,x'_1,x'_3,z_4) = (16,5,8,4,4,2,2,2,1)$ and\\ 
    $\bullet$
    $\bT_3(z_5,x_5,z'_5,z'_1,x'_5,x'_1,x'_3,x'_4,z_3) = (16,5,8,4,4,2,2,1,1)$
\end{quote}
giving rise to the inequality $|C| \le \frac{1}{5}(\bT_1+\bT_2+\bT_3) \le 133/5$ whenever $C$ is $z_1$-unsolvable.
Hence $\pi(B_2,z_1)\le 27$.

For $B_2$, we define three $z_2$-strategies $\bT_1$, $\bT_2$, and $\bT_3$ by 
\vspace{-0.04in}
\begin{quote}
    $\bullet$ 
    $\bT_1(z'_2,z'_1,x'_2,z'_5,x'_1,x'_4,x'_5,x'_3) = (16,8,8,4,4,4,4,2)$, \\
    $\bullet$
    $\bT_2(z_1,z_5,x_1,x_3,z_3,x'_3) = (16,4,8,4,2,1)$ and\\ 
    $\bullet$
    $\bT_3(x_2,x_5,x_4,z_4,z_3,x'_3) = (16,4,8,4,2,1)$
\end{quote}
giving rise to the inequality $|C| \le \frac{1}{4}(\bT_1+\bT_2+\bT_3) \le 120/4$ whenever $C$ is $z_2$-unsolvable.
Hence $\pi(B_2,z_2)\le 31$.

For $B_2$, we define three $z_3$-strategies $\bT_1$, $\bT_2$, and $\bT_3$ by 
\vspace{-0.04in}
\begin{quote}
    $\bullet$ 
    $\bT_1(z_4,x_4,x'_4,x_2,x'_2,z_2,z'_2) = (32,16,16,8,8,4,4)$, \\
    $\bullet$
    $\bT_2(x'_3,x'_1,x'_5,z'_1,z'_5,z'_2) = (32,16,16,8,8,4)$ and\\ 
    $\bullet$
    $\bT_3(x_3,x_1,x_5,z_1,z_5,z_2) = (32,16,16,8,8,4)$
\end{quote}
giving rise to the inequality $|C| \le \frac{1}{8}(\bT_1+\bT_2+\bT_3) \le 133/5$ whenever $C$ is $z_3$-unsolvable.
Hence $\pi(B_2,z_3)\le 33$.

\end{document}